\documentclass[final]{article}

\usepackage{amsmath,amssymb,amsthm}
\newtheorem{theorem}{Theorem}
\newtheorem{lemma}{Lemma}

\usepackage{amsfonts}
\usepackage{epsfig}

\newcommand{\Laplacian}{\mathrm{\Delta}}
\newcommand{\ba}{\begin{array}}
\newcommand{\ea}{\end{array}}
\newcommand{\be}{\begin{equation}}
\newcommand{\ee}{\end{equation}}
\newcommand{\Id}{\mathbb{I}}
\newcommand{\st}{\sigma_{\tau}}
\newcommand{\FE}{F_E}
\newcommand{\TV}{\mathrm{TV}}

\title{High order relaxation schemes for non linear degenerate
        diffusion problems} 

\author{Fausto Cavalli\footnotemark[2]{~}\footnotemark[4] 
   \and Giovanni Naldi \footnotemark[2]{~}\footnotemark[4] 
   \and Gabriella Puppo\footnotemark[3]{~}\footnotemark[4] 
   \and Matteo Semplice\footnotemark[2]{~}\footnotemark[4] }

\begin{document}

\maketitle
\renewcommand{\thefootnote}{\fnsymbol{footnote}}
\footnotetext[2]{Dipartimento di Matematica, Universit\`a di Milano,
Via Saldini 50, I-20100 Milano, Italy.  
}
\footnotetext[3]{Dipartimento di Matematica, Politecnico di Torino, 
 Corso Duca degli Abruzzi 24, 10129 Torino, Italy.
{\rm gabriella.puppo@polito.it}}
\footnotetext[4]{This work was partially supported by the MIUR/PRIN2005 project
``Modellistica numerica per il calcolo scientifico ed applicazioni
avanzate''.}
\renewcommand{\thefootnote}{\arabic{footnote}}

\begin{abstract}
Several relaxation approximations to partial differential equations
have been recently proposed. Examples include conservation laws,
Hamilton-Jacobi equations, convection-diffusion problems, gas dynamics
problems. The present paper focuses onto diffusive relaxation schemes
for the numerical approximation of nonlinear parabolic
equations. These schemes are based on a suitable semilinear hyperbolic
system with relaxation terms. High order methods are obtained by
coupling ENO and WENO schemes for space discretization with IMEX
schemes for time integration.  Error estimates and convergence
analysis are developed for semidiscrete schemes with numerical
analysis for fully discrete relaxed schemes.  Various numerical
results in one and two dimensions illustrate the high accuracy and good
properties of the proposed numerical schemes, also in the degenerate
case. These schemes can be easily implemented on parallel computers
and applied to more general systems of nonlinear parabolic equations
in two- and three-dimensional cases.
\end{abstract}

\section{Introduction}                  \label{section:introduction}
\setcounter{equation}{0}
\setcounter{figure}{0}
\setcounter{table}{0}

Relaxation approximations to nonlinear partial differential equations
have been introduced \cite{AN00,ANT04,JPT98,JX95,LT97,NP00,NPT02} on
the basis of the replacement of the equations with a suitable
semilinear hyperbolic system with stiff relaxation terms. The idea can
be explained by considering the case of a scalar conservation law
\[
\frac{\partial u}{\partial t} + \frac{\partial f(u)}{\partial x} =0,
\]
for which Jin and Xin proposed in \cite{JX95} the following relaxation
approximation
\[
\left\{
\begin{array}{ll}
\displaystyle
\frac{\partial u}{\partial t} + \frac{\partial v}{\partial x} = 0\\
\\
\displaystyle
\frac{\partial v}{\partial t} +  a  \frac{\partial u}{\partial t}= 
    -\frac{1}{\epsilon}\left( v-f(u)\right)
\end{array}
\right.
\]
where $\epsilon$ (relaxation time) is a small positive parameter and
$a$ is a positive constant satisfying  
\[
-\sqrt{a} \leq f^\prime (u) \leq \sqrt{a}
\]
for all $u$. For small values of $\epsilon$ the previous system gives
the following first order approximation  of the original conservation
law 
\[ 
\frac{\partial u}{\partial t} + \frac{\partial f(u)}{\partial x} = \epsilon \frac{\partial}{\partial x} \left(
(a-f^\prime (u)^2) \frac{\partial u}{\partial x}\right).
\]
The generalization of the above model to multidimensional systems of
conservation laws can be done in a natural way by adding more rate
equations. One would expect that appropriate numerical schemes for the
relaxation system yield accurate numerical approximations to the
original equation or system when the relaxation rate $\epsilon$ is
sufficiently small. Numerically, the main advantage of solving the
relaxation model over the original conservation law lies in the simple
linear structure of characteristic fields and the localized lower
order term. In particular, the semilinear nature of the relaxation
system gives a new way to develop numerical schemes that are simple,
general and Riemann solver free \cite{Jin95,JL96,JX95}.  Several other
relaxation approximation have been introduced recently. For example we
mentioned here the work of Coquel and Perthame \cite{CP98} for real
gas computation, the relaxation schemes of Jin et al. \cite{JKX99} for
curvature-dependent front propagation, the relaxation approximation in
the rapid granular flow \cite{JS01}, the relaxation approximation and
relaxation schemes for diffusion and convection-diffusion problems
\cite{JPT98,LT97,NP00,NPT02}. Moreover, there are strong and interesting
links between the relaxation approximation and the kinetic approach to
nonlinear transport equations, based upon analogies with the passage
from the Boltzmann equation to fluid mechanics (see for example
\cite{AN00,ANT04,BGN00}).

The aim of this work is to analyze from both a theoretical and
computational point of view relaxation schemes that approximate the
following nonlinear degenerate parabolic problem
\begin{equation}
\label{eq:degparab}
\frac{\partial u}{\partial t} = D \Laplacian (p(u)),~~~x\in \mathbb{R}^d,~t>0, 
\end{equation}
with initial data $u(x,0)=u_0(x)\in L^1(\mathbb{R}^d)$ and $D>0$ is a
diffusivity coefficient. As usual, we will assume $p: \mathbb{R}
\rightarrow \mathbb{R}$ to be  non-decreasing and
Lipschitz continuous \cite{Vaz92}. The equation is degenerate if
$p(0)=0$. We will also consider the same parabolic problem in a
bounded domain $\Omega \subset \mathbb{R}^d$ with Neumann boundary
conditions.  In the case $p(u)=u^m$, $m>1$ the previous equation is
the \textit{porous media equation} which describes the flow of a gas
through a porous interface according to some constitutive relation
like Darcy's law in order to link the velocity of the gas and its
pressure.  In this case the diffusion coefficient $mu^{m-1}$ vanishes
at the points where $u\equiv 0$ and the governing parabolic equation
degenerates there. The set of such points is called
interface. Moreover the porous media equation can exhibit a finite
speed of propagation for compactly supported initial data
\cite{Aro70}.  The influence of the degenerate diffusion terms
make the dynamics of the interfaces difficult to study from both the
theoretical and the numerical point of view.  Another interesting case
corresponds to $0<m<1$ and it is referred to as the \textit{fast
diffusion equation} which appears, for example, in curvature-driven
evolution and avalanches in sandpiles.
In general the numerical analysis of
equation \eqref{eq:degparab} is difficult for at least two reasons:
the appearance of singularities for compactly supported solutions and the
growth of the size of the support as time increases (\textit{retention
property}).

From the numerical viewpoint, an usual technique to approximate
\eqref{eq:degparab} involves implicit discretization in time: it
requires, at each time step, the discretization of a nonlinear
elliptic problem. However, when dealing with nonlinear problems one
generally tries to linearize them in order to take advantage of
efficient linear solvers. Linear approximation schemes based on the
so-called nonlinear Chernoff's formula with a suitable relaxation
parameter have been studied for example in
\cite{BBR79,NSV00,NV88,MNV87} where also some energy error estimates
have been investigated. Other linear approximation schemes have been
introduced by J{\"a}ger, Ka{\v{c}}ur and Handlovi{\v{c}}ov{\'a}
\cite{JK1,Kac1}. More recently, different approaches based on
kinetic schemes for degenerate parabolic systems have been considered and
analyzed by Aregba-Driollet, Natalini and Tang in
\cite{ANT04}. Other approaches were investigated in the work of
Karlsen et al. \cite{EK99} based on a suitable splitting technique with
applications to more general hyperbolic-parabolic convection-diffusion
equations. Finally, a new scheme based on the maximum principle and on
a perturbation and regularization approach was proposed by Pop and
Yong in \cite{PY02}.

Our approach is inspired by relaxation schemes where the nonlinearity
inside the equation is replaced by a semilinearity. This reduction is
carried out in order to obtain numerical schemes that are easy to implement,
also for parallel computing, even in the multidimensional case and for
more general and complex problems, like oil recovery problems \cite{EK00}.
Moreover with our approach it is
possible to improve numerical schemes by using high order methods and,
in principle, different numerical approaches (finite volume, finite
differences,\ldots). In particular, in this paper, in order to obtain
high order methods, we couple ENO and WENO schemes for space
discretization and IMEX schemes for time advancement.
High order schemes may not reach their order of convergence due to the
loss of regularity of the solution during the evolution. However they
are interesting nevertheless for error reduction when the number of
grid points is fixed or until discontinuities develop (both cases
arise for example in nonlinear filtering in image analysis \cite{Wei98}).

The paper is organized as follows. Section \ref{section:relaxation} is
devoted to the introduction of our relaxation schemes. The stability
and error estimates of the semidiscrete scheme is provided in Section
\ref{section:semidiscrete}.  In Section \ref{section:fullydiscrete} we
consider the fully discrete relaxed scheme with nonlinear stability
analysis, the numerical treatment of boundary conditions and the
extension to the multi-dimensional case. Finally, the implementation
of the method as well as the results of several numerical experiments
are discussed in section \ref{section:numerical}. 
 
\section{Relaxation approximation of nonlinear diffusion}                  
\label{section:relaxation}
\setcounter{equation}{0}
\setcounter{figure}{0}
\setcounter{table}{0}

The main purpose of this work is to approximate solutions of a nonlinear,
possibly degenerate, parabolic equation of the form \eqref{eq:degparab},
This framework is so general that it includes the
porous medium equation $p(u)=u^m$ with $m >1$ (\cite{Vaz92} and
references therein), 
non linear image processing \cite{Wei98}, as well as a wide class of mildly
nonlinear parabolic equations \cite{Fri60}.

The schemes proposed in the present work are based on the same idea at
the basis of the well-known relaxation schemes for hyperbolic
conservation laws \cite{JX95}.  In the case of the nonlinear diffusion
operator, an additional variable $\vec{v}(x,t)\in\mathbb{R}^d$ and a
positive parameter $\varepsilon$ are introduced and the following
relaxation system is obtained:
\begin{equation}
\label{sysrel1}
\left\{
\begin{array}{ll}
\displaystyle
\frac{\partial u}{\partial t} + \mathrm{div}(\vec{v}) = 0\\
\\
\displaystyle
\frac{\partial \vec{v}}{\partial t} + \frac{D}{\varepsilon} \nabla p(u) = 
    -\frac1\varepsilon \vec{v}
\end{array}
\right.
\end{equation}
Formally, in the small relaxation limit, $\varepsilon \rightarrow
0^+$, system \eqref{sysrel1} approximates to leading order
equation \eqref{eq:degparab}.
In order to have bounded characteristic velocities and to avoid a singular differential operator as 
$\varepsilon \rightarrow 0^+$, a suitable parameter $\varphi$ is introduced and
\eqref{sysrel1} can be rewritten as: 
\begin{equation}
\left\{
\begin{array}{ll}
\displaystyle
\frac{\partial u}{\partial t} + \mathrm{div}(\vec{v}) = 0\\
\\
\displaystyle
\frac{\partial \vec{v}}{\partial t} + \varphi^2 \nabla p(u) = 
    -\frac1\varepsilon \vec{v} + \left(\varphi^2-
    \frac{D}{\varepsilon}\right) \nabla p(u)
\end{array}
\right.
\end{equation}
Finally we remove the non linear term from the convective part, as in
standard relaxation schemes, introducing a variable
$w(x,t)\in\mathbb{R}$ and rewriting the system as:
\begin{equation} \label{3eq}
\left\{
\begin{array}{ll}
\displaystyle
\frac{\partial u}{\partial t} + \mathrm{div}(\vec{v}) = 0 \\
\\
\displaystyle \frac{\partial \vec{v}}{\partial t} + \varphi^2 \nabla w = 
    -\frac1\varepsilon \vec{v}
    +\left(\varphi^2- \frac{D}{\varepsilon}\right) \nabla w \\
\\
\displaystyle
\frac{\partial w}{\partial t} + \mathrm{div}(\vec{v}) =
    -\frac1\varepsilon (w-p(u))
\end{array} 
\right.
\end{equation}
Formally, as $\varepsilon \rightarrow 0^+$, $w \rightarrow p(u)$, $v
\rightarrow -\nabla p(u)$ and the original equation is recovered. 

In the previous system the parameter $\varepsilon$ has physical
dimensions of time and represents the so-called relaxation
time. Furthermore, $w$ has the same dimensions as $u$, while each
component of $\vec{v}$ has the dimension of $u$ times a velocity;
finally $\varphi$ is a velocity.  The inverse of $\varepsilon$ gives the
rate at which $v$ decays onto $-\nabla p(u)$ in the evolution of the
variable $\vec{v}$ governed by the stiff second equation of
\eqref{3eq}.

Equations \eqref{3eq} form a semilinear hyperbolic system with a stiff
source term. The characteristic velocities of the hyperbolic part are
given by $0,\pm\varphi$. The parameter $\varphi$ allows to move the stiff
terms $\frac{D}{\varepsilon}\nabla p(u)$ to the right hand side,
without losing the hyperbolicity of the system.

We point out that degenerate parabolic equations often model physical
situations with free boundaries or discontinuities: we expect that
schemes for hyperbolic systems will be able to reproduce faithfully
these details of the solution.  One of the main properties of
\eqref{3eq} consists in the semilinearity of the system, that is all
the nonlinearities are in the (stiff) source terms, while the
differential operator is linear. Hence, the solution of the convective
part requires neither Riemann solvers nor the computation of the
characteristic structure at each time step, since the eigenstructure
of the system is constant in time.  Moreover, the relaxation
approximation does not exploit the form of the nonlinear function $p$
and hence it gives rise to a numerical scheme that, to a large extent,
is independent of it, resulting in a very versatile tool.

We also anticipate here that, in the relaxed case
(i.e. $\varepsilon=0$), the stiff source terms can be integrated
solving a system that is already in triangular form and then it does
not require iterative solvers.
\section{The semidiscrete scheme}                  
\label{section:semidiscrete}
\setcounter{equation}{0}
\setcounter{figure}{0}
\setcounter{table}{0}

System \eqref{3eq} can be written in the form:
\be \label{eq:statement}
z_t + \mathrm{div} f(z) = \frac{1}{\varepsilon}g(z),  
\ee
where 
\be
 z=\left(\ba{c}u\\v\\w\ea\right)
 \qquad
 f(z)=\left[\ba{c}v^T\\\Phi^2w\\v^T\ea\right]
 \qquad
 g(z)=\left(\ba{c}0\\ -v +(\varphi^2 \varepsilon -D) \nabla w\\ p(u)-w
\ea\right)
\ee
and $\Phi^2$ is the $d\times d$ identity matrix times the scalar 
$\varphi^2$. We start discretizing the
system in time using, for simplicity, a uniform time step $\Delta t$. Let $z^n(x) =
z(x,t^n)$, with $t^n = n\Delta t$. Since equation \eqref{eq:statement}
involves both stiff and non-stiff terms, it is a natural idea to
employ different time-discretization strategies for each of them, as
in \cite{ARS97,PR05}. In this work we integrate \eqref{eq:statement}
with a Runge-Kutta IMEX scheme \cite{PR05}, obtaining the following
semidiscrete formulation
\be \label{eq:balance:imex}
z^{n+1} = z^n - \Delta t \sum_{i=1}^{\nu}
     \tilde{b}_i  \mathrm{div} f(z^{(i)}) 
    + \frac{\Delta t}{\varepsilon} \sum_{i=1}^{\nu} b_i
     g(z^{(i)}),  
\ee
where the $z^{(i)}$'s are the stage values of the Runge-Kutta scheme
which are given by
\be
z^{(i)} = z^n -\Delta t \sum_{k=1}^{i-1}\tilde{a}_{i,k}
   \mathrm{div} f(z^{(k)}) + \frac{\Delta t}{\varepsilon}
         \sum_{k=1}^{i} a_{i,k} g(z^{(k)}),
\label{eq:balance:intermediate} \ee 
where $\tilde{b}_i$, $\tilde{a}_{ij}$ and $b_i$, $a_{ij}$  denote the coefficients of the
explicit and implicit RK schemes, respectively. We assume that the
implicit scheme is of diagonally implicit type. To find the $z^{(i)}$'s it is
necessary in principle to solve a non linear system of equations which
however can be easily decoupled. The system for the first stage
$z^{(1)}$ at time $t^n$ is:
\be \label{eq:stage:1}
\left(\ba{c} u^{(1)} \\ v^{(1)} \\ w^{(1)}\ea \right)=\left(\ba{c} u^{n} \\ v^{n} \\ w^{n}\ea \right)+\frac{\Delta t}{\varepsilon}a_{11}\left(\ba{c} 0 \\ -v^{(1)}+(\varphi^2\varepsilon-D)\nabla w^{(1)} \\ p(u^{(1)})-w^{(1)}\ea \right).
\ee
The first equation yields $u^{(1)}=u^{n}$, substituting in the third
equation we immediately find $w^{(1)}$ and finally, substituting $w^{(1)}$
in the second equation, we compute  $v^{(1)}.$ In other words the system can
be written in triangular form.  For the following stage values, grouping
the already computed terms in the vector $B^{(i)}$ given by
\be \label{eq:stage:expl}
B^{(i)} = z^n -\Delta t \sum_{k=1}^{i-1}\tilde{a}_{i,k}
   \mathrm{div} f(z^{(k)}) + \frac{\Delta t}{\varepsilon}
         \sum_{k=1}^{i-1} a_{i,k} g(z^{(k)}),
\ee
then the new stage values are given by
\be \label{eq:stage:i}
\left(\ba{c} u^{(i)} \\ v^{(i)} \\ w^{(i)}\ea \right)=B^{(i)}+\frac{\Delta t}{\varepsilon}a_{ii}\left(\ba{c} 0 \\ -v^{(i)}+(\varphi^2\varepsilon-D)\nabla w^{(i)} \\ p(u^{(i)})-w^{(i)}\ea \right),
\ee
which is again a triangular system.  In the numerical tests, we will
apply IMEX schemes of order 1, 2 and 3. 

Following \cite{JX95} we set $\varepsilon=0$ thus obtaining the so
called \textit{relaxed scheme}. The computation of the first stage
reduces to
\be\label{firststageRK}
\ba{l} u^{(1)}=u^n \\ w^{(1)}=p(u^{(1)}) \\ v^{(1)}=-D\nabla w^{(1)}\ea.
\ee
For the following stages the first equation is
\be\label{originalRK}
 u^{(i)}=u^n -\Delta t \sum_{k=1}^{i-1}\tilde{a}_{i,k}
   \mathrm{div} v^{(k)}.
\ee
In the other equations the convective terms are dominated by the
source terms and thus $v^{(i)}$ and $w^{(i)}$ are given by
\be\label{relaxedRK}
\ba{l}
 v^{(i)}=-D\nabla w^{(i)}, \\
 w^{(i)}=p(u^{(i)}).
\ea
\ee

We see that only the explicit part of the Runge-Kutta method is involved in the updating
of the solution. Then, in the relaxed schemes we use only the explicit part of the tableaux. In particular we consider  second and third order SSRK schemes
\cite{GST01}, namely
\begin{center}
\begin{tabular}{c@{\hspace{1.5cm}}c@{\hspace{1.5cm}}c}
{\tt IMEX1} ($1^{\mbox{st}}$ order) 
&{\tt IMEX2} ($2^{\mbox{nd}}$ order)  
&{\tt IMEX3} ($3^{\mbox{rd}}$ order)\\[3mm]
\raisebox{-12pt}{
$\begin{array}{c|c}
        & 0\\
\hline\\[-3mm]  & 1
\end{array}
$}
&
\raisebox{-6pt}{
$
\begin{array}{c|cc}
  & 0& 0\\
  & 1& 0\\
\hline \\[-3mm]
   & \frac12& \frac12
\end{array}
$}
&
$\begin{array}{c|cccc}
 &0& 0& 0\\
 &1 &0 &0\\
 &\frac14& \frac14& 0\\
 \hline \\[-3mm]
   & \frac16 &\frac16 & \frac23
\end{array}
$
\end{tabular}
\end{center}

\subsection{Convergence of the semidiscrete relaxed scheme}
The aim of this section is to show the $L^1$ convergence of the
solution of the semidiscrete in time relaxed scheme defined by
equations \eqref{firststageRK},\eqref{originalRK} and
\eqref{relaxedRK}. We will extend the theorem proved in \cite{BBR79},
where only the case of forward Euler timestepping was
considered. In this section, for the sake of simplicity, we set
$D=1$.

Eliminating $v$ from \eqref{firststageRK} and \eqref{originalRK} using
\eqref{relaxedRK}, we rewrite the relaxed scheme as
\be
\ba{l} u^{(1)}=u^n \\ w^{(1)}=p(u^n)\ea
\ee
for the first stage, and 
\be
\ba{l}  u^{(i)}=u^n +\Delta t \sum_{k=1}^{i-1}\tilde{a}_{i,k}
        \mathrm{\Delta} w^{(k)}\\
        w^{(i)}=p(u^{(i)}),
\ea\ee
for subsequent stages.  
We recall that a Runge-Kutta scheme for the ordinary differential equation \(y'=R(y)\) can
also be written in the form \cite{GST01}
\be\label{imex2bis}
\ba{l}
  y^{(1)}=y^n\\
  y^{(i)}=\sum_{k=1}^{i-1}\alpha_{ik}\left(y^{(k)}+\Delta t\frac{\beta_{ik}}{\alpha_{ik}}R(y^{(k)})\right)  \qquad i=2,\dots, \nu,
\ea\ee
where \(y^{n+1}=y^{(\nu)}\).
For consistency,
\(\sum_{k=1}^{i-1}\alpha_{ik}=1\) for every
\(i=1,\ldots,\nu\). Moreover we assumed that \(\alpha_{ik}\geq0\),
\(\beta_{ik}\geq0\) and that \(\alpha_{ik}=0\) implies
\(\beta_{ik}=0\). 
Under these assumptions, each stage value \(y^{(i)}\) can be written as a convex
combination of forward Euler steps. This remark allows us to study the
convergence of the Runge-Kutta scheme in terms of the convergence of
the explicit forward Euler scheme applied to the non-linear diffusion
problem.

This latter was studied in \cite{BBR79} via a nonlinear semigroup
argument.  In the following we review the approach of \cite{BBR79} and
next we extend the proof to the case of a $\nu$-stages explicit
Runge-Kutta scheme.

\subsubsection{The forward Euler case}
We wish to solve the evolution equation
\be\label{eveq}
	\frac{du}{dt}+Lp(u)=0  \qquad u(\cdot,t=0)=u_0,
\ee
on the domain $\Omega$, where $L=-\Delta$ and
$p:\mathbb{R}\rightarrow\mathbb{R}$ is a non decreasing locally
Lipschitz function such that $p(0)=0$.
Under these hypotheses, the nonlinear operator $Au=Lp(u)$ with domain
$D(A)=\{u\in L^1(\Omega): p(u)\in D(L)\}$ is m-accretive in
$L^1(\Omega)$, that is $\forall\varphi\in\mathrm{L}^1(\Omega)$ and
$\forall\lambda>0$ there exists a unique solution $u\in D(A)$ such
that \(u+\lambda Lp(u)=\varphi\) and the application defined by
\(\varphi\mapsto u\) is a contraction\cite{CL71}.

Moreover $D(A)$ is dense in $L^1(\Omega)$, so it
follows that
\be
  S_A(t)u_0=\lim_{m\rightarrow\infty}\left(\Id+\frac{t}{m}A\right)^{-m}u_0
\ee
is a contraction semigroup on $L^1(\Omega)$ and $S_A(t)u_0$ is the
generalized solution of \eqref{eveq} in the sense of Crandall-Liggett
\cite{CL71}. 
Let $S(t)$ be the linear contraction semigroup generated by $-L$, that
is $u(t)=S(t)u_0$ is the solution of the initial value problem
$u_t=-L(u)$ and $u(\cdot,t=0)=u_0$.
The algorithm proposed in \cite{BBR79} is
\be\label{Brezis1}
\frac{u^{n+1}-u^n}{\tau} + 
\left[\frac{\Id-S(\st)}{\st}\right]p(u^n)=0,
\ee
where $\tau$ is the timestep and $\st\downarrow0$.
This can be written as
\be \label{Brezis:onestep}
u^{n+1}=\FE(\tau)u^n \quad\mbox{ where }
\FE(\tau)\varphi=\varphi+\frac\tau\st\left[S(\st)-\Id\right]p(\varphi).
\ee
Hence
\be
u^n=(\FE(\tau))^n u_0.
\ee
The proof in \cite{BBR79} is based on the following argument. Note that formally
\(S(\st)\varphi\sim e^{-\st L}\varphi\). Let $t=\tau n$
\be \label{Brezis:argument}
\ba{ll}
  u(t) &= \left[\Id+\frac{t}{n\st}\left(S(\st)-\Id\right)\circ
         p\right]^n u_0\\
       &= \left[\Id+\frac{t}{n\st}\left(e^{-\st L} -\Id\right)\circ
         p\right]^n u_0 \qquad\mbox{if $\st\rightarrow0$}\\
       &= \left[\Id-\frac{t}{n}L\circ p\right]^n u_0\\
       &\rightarrow S_A(u_0) \qquad\mbox{when $n\rightarrow\infty$}
\ea
\ee
The convergence proof requires that $\mu\frac\tau\st\leq1$ where $\mu$
is the Lipschitz constant of $p(u)$. We point out that $\st$ is linked
to the spatial approximation of the operator $L$ and in our scheme
this requirement is reflected in the stability condition of the fully
discrete scheme (see Section \ref{section:fullydiscrete})

\subsubsection{Runge-Kutta schemes}

Now we are going to describe the case of a $\nu$-stages Runge-Kutta
scheme, proving its convergence.

Let $t>0$ and $\tau=t/n$ with $n\geq1$; let $\st:(0,\infty)\rightarrow
(0,\infty)$ be a function such that $\lim_{\tau \rightarrow0} \st=0$.

\be
\ba{l}
  u^{(1)}=u^n,\\
  u^{(i)}=\sum_{k=1}^{i-1}\alpha_{ik}\left[
                u^{(k)}+\tau\frac{\beta_{ik}}{\alpha_{ik}}
                A(u^{(k)}) 
                \right]
\qquad i=2,\ldots,\nu
\ea\ee
and proceeding as in (\ref{Brezis:argument}), this becomes
\be\label{Brezis:RK}
\ba{l} 
  u^{(1)}=u^n,\\
  u^{(i)}=\sum_{k=1}^{i-1}\alpha_{ik}\left[
                u^{(k)}+\tau\frac{\beta_{ik}}{\alpha_{ik}}
                (S(\st)-\Id)\circ{p(u^{(k)})}   
                \right] 
          \qquad i=2,\ldots,\nu\\
  u^{n+1}=u^{(\nu)}
\ea\ee

We now extend \eqref{Brezis:onestep} to the Runge-Kutta scheme defined
by equation \eqref{Brezis:RK}. Define, for $\phi\in L^1(\Omega)$,
\be \label{Brezis:F}
\ba{l}
   F^{(1)}(\tau)\phi=\phi, \\
   F^{(i)}(\tau)\phi=\sum^{i-1}_{k=1}\alpha_{ik}F^{(k)}(\tau)\phi
       +\frac{\tau\beta_{ik}}{\st}\left[S(\st)-\Id\right]p(F^{(k)}(\tau)\phi), \\
   F(\tau)\phi=F^{(\nu)}(\tau)\phi
\ea\ee
and therefore
\be
	u^n(t)=\left[ F(\tau)\right]^n u_0.
\ee

Let $u(t)$ be the generalized solution of \eqref{eveq}. The following
theorem proves the convergence of the semidiscrete solution to $u(t)$.
\begin{theorem} \label{teorema}
	Assume $u^0\in L^{\infty}(\Omega),$ and $\|u^0\|_{\infty}=M$;
	let $p$ be a non-decreasing Lipschitz continuous function on
	$[-M,M]$ with Lipschitz constant $\mu$. Assume that the
	following conditions hold 
\begin{equation}\label{stab}
	\left\{\begin{array}{l} \displaystyle \alpha_{ik}\geq0, \\
	\displaystyle \beta_{ik}\geq0, \\ \displaystyle
	\alpha_{ik}=0\Rightarrow\beta_{ik}=0, \\ \displaystyle
	\sum^{i-1}_{k=1}\alpha_{ik}=1 \mbox{ (consistency)},\\[5mm]
	\displaystyle 
        \frac{\mu\tau}{\st}\leq \min\frac{\alpha_{ik}}{\beta_{ik}}, 
         \quad \mbox{for }\tau>0,\alpha_{ik}\neq 0 \
         \quad \mbox{ (stability)},
	\end{array}\right.  
\end{equation} 
then
	$\lim_{n\rightarrow\infty} u^n(t)=u(t)$ in $L^1$. Moreover the
	convergence is uniform for t in any given bounded interval.
\end{theorem}

The proof follows the steps of \cite{BBR79}: first we show that $u^n$
verifies a maximum principle (Lemma \ref{lemma1}) and that $F$ is a
contraction (Lemma \ref{lemma2}) and finally we apply the non linear
Chernoff formula\cite{BP72}.

\begin{lemma} \label{lemma1}
 If \eqref{stab} is verified, then $-M\leq u^n\leq M \quad \forall n.$
\end{lemma}
\begin{proof}
We argue by induction on $n$: we assume that $-M\leq u^n\leq M$ and we
show that $-M\leq u^{n+1}\leq M$. Let
\begin{equation}\label{inter}
u^{(i)}=F^{(i)}(\tau)u^{n}
\end{equation}
Since $u^{n+1}=u^{(\nu)}$, it suffices to prove that $-M\leq
u^{(i)}\leq M$ for $i=1,...,\nu$. We prove this by induction on $i$. When
$i=1$, the statement is true thanks to the induction hypothesis on $n$ and
being $F^{(1)}=\Id$.
Let's assume that $-M\leq u^{(i-1)}\leq M$ holds; we are going to show that
\begin{equation}\label{ii}
-M\leq u^{(i)}=F^{(i)}(\tau)u^{n}\leq M.
\end{equation}

The function \(s\mapsto \alpha_{ik}s-\frac{\tau\beta_{ik}}{\st}p(s)\)
is non decreasing thanks to \eqref{stab} and the hypotheses on the
function $p$. By the induction hypothesis on $i$, we have that for $k=1, ...,i-1$
\begin{equation} \label{induction:i}
 -\alpha_{ik}M-\frac{\tau\beta_{ik}}{\st}p(-M)\leq
  \alpha_{ik}u^{(k)}-\frac{\tau\beta_{ik}}{\st}p(u^{(k)})\leq
  \alpha_{ik}M-\frac{\tau\beta_{ik}}{\st}p(M).
\end{equation}
Using again the induction hypothesis on $i$, recalling that $p$ is
non-decreasing, since $S$ is a
contraction in $\mathrm{L}^\infty$ \cite{BBR79} and \(p(-M)\leq p(u^{(k)}) \leq p(M)\),
\be
  p(-M)\leq S\left(p(u^{(k)})\right) \leq p(M)
\ee
Multiplying the last equation by $\frac{\tau\beta_{ik}}{\st}$ and
summing it to equation \eqref{induction:i}, we get
\begin{equation} \label{induction:ik}
 -\alpha_{ik}M\leq
  \alpha_{ik}u^{(k)}+\frac{\tau\beta_{ik}}{\st}(S-\Id)p(u^{(k)})\leq
  \alpha_{ik}M,
\qquad k=1,\ldots,i-1
\end{equation}
Summing for $k=1,...,i-1$ and using the consistency relation of \eqref{stab}:
\begin{equation}\label{ai}
   -M
  \leq \sum^{i-1}_{k=1}\alpha_{ik}u^{(k)}+\frac{\tau\beta_{ik}}{\st}(S-\Id)p(u^{(k)})
  \leq M
\end{equation}
In particular this is valid when $i=\nu$, proving that \(-M\leq
u^{(n+1)}\leq M\).
\end{proof}

Now we can replace $p$ by $\overline{p}$, where $\overline{p}=p$ in
$-M\leq x\leq M,$ $\overline{p}=p(M)$ for $x\geq M$ and
$\overline{p}=p(-M)$ for $x\leq -M:$ the algorithm is the same and in
what follows we can assume that $p$ is Lipschitz continuous with
constant $\mu$ on all $\mathbb{R}$.

\begin{lemma}\label{lemma2}
 If the hypotheses of Theorem \ref{teorema} hold, then $F(\tau)$
 is a contraction on $L^1(\Omega),$ i.e.  
 \begin{equation}\label{lem2}
   \|F(\tau)\phi-F(\tau)\psi\|_1\leq \|\phi-\psi\|_1 \qquad \forall
   \psi,\phi\in L^1
 \end{equation}
\end{lemma}
\begin{proof}
We start showing that the result holds for a single forward Euler
step. Recalling the definition of \(\FE\) from \eqref{Brezis:onestep}
\be\ba{ll}
  \|\FE(\tau)\phi-\FE(\tau)\psi\|_1
  &\leq
   \frac{\tau}{\st}\left\|S(\st)[p(\phi)-p(\psi)]\right\|_1
   +\left\|(\phi-\psi)-\frac{\tau}{\st}[p(\phi)-p(\psi)]\right\|_1 \\
  &\leq
   \frac{\tau}{\st}\left\|p(\phi)-p(\psi)\right\|_1
   +\left\| \left(\phi-\frac{\tau}{\st}p(\phi)\right)
           -\left(\psi-\frac{\tau}{\st}p(\psi)\right)
    \right\|_1\\
  &=\left\|\phi-\psi\right\|_1
\ea\ee
where we used the contractivity of $S$.
The last equality relies on the fact that $p$ and the function
$x\mapsto{}x-\frac{\tau}{\st}p(x)$ are non-decreasing, which in turn is
guaranteed by the stability condition, that in this case
reduces to \(\mu\tau/\st\leq1\) \cite{BBR79}.

In the general case we have:
\be\label{uno}
\ba{ll}
\|F^{(i)}(\tau)\phi-F^{(i)}(\tau)\psi\|_1
 &\leq \sum_{k=1}^{i-1}\alpha_{ik}
    \left\|\FE\left(\frac{\tau\beta_{ik}}{\alpha_{ik}}\right)F^{(k)}(\tau)\phi
           -\FE\left(\frac{\tau\beta_{ik}}{\alpha_{ik}}\right)F^{(k)}(\tau)\psi
    \right\|_1\\
 &\leq \sum_{k=1}^{i-1}\alpha_{ik}
    \left\|F^{(k)}(\tau)\phi-F^{(k)}(\tau)\psi\right\|_1\\
 &\leq \left\|\phi-\psi\right\|_1
\ea\ee
In the second inequality we used the contractivity of $\FE$ and the
stability condition, while in the third one we apply an induction
argument on the contractivity of $F^{(k)}$, the positivity
constraint on $\alpha_{ik}$ and $\beta_{ik}$, as well as the
consistency condition $\sum_k\alpha_{ik}=1$. Setting $i=\nu$ yields the result.
\end{proof}

\begin{proof}[Proof of Theorem \ref{teorema}]
Let \(\psi_{\tau}\) and \(\psi\) be respectively
\be
	\psi_{\tau}=\left(I+\frac{\lambda}{\tau}(I-F(\tau))\right)^{-1}\phi
\qquad\mbox{and}\qquad	
	\psi=\left(I+\lambda A\right)^{-1}\phi. 
\ee
The function $\psi$ exists since the operator $A$ is m-accretive,
whereas the existence of the function $\psi_\tau$ is guaranteed by the
following fixed-point argument. Let
\[ G(y) = \frac1{1+\eta}\phi+\frac{\eta}{\eta+1}F(\tau)y \]
where $\phi\in\mathrm{L}^1$, $y\in\overline{D(A)}$ and
$\eta\geq0$. We have, 
\[ \|G(y)-G(x)\|=\frac{\eta}{\eta+1}\|F(\tau)y-F(\tau)x\|
   \leq \frac{\eta}{\eta+1}\|y-x\| 
\] 
since $F$ is a contraction, as proved in Lemma \ref{lemma2}.
Thus $G$ is also a contraction and therefore it possesses a unique
fixed point which coincides with  $\psi_\tau$.

We want to show that 
\[ \psi_{\tau} \rightarrow \psi 
   \qquad
   \mbox{in } \mathrm{L}^1
\]
as \(\tau\rightarrow0\) for each fixed $\lambda>0$.
 Let
\[ \phi_\tau = \psi+\frac\lambda\tau (\Id - F(\tau))\psi. \]
We want to estimate \(\psi_\tau-\psi\) in terms of
\(\phi_\tau-\phi\).
\[
  \phi_\tau-\phi = (1+\frac\lambda\tau)(\psi-\psi_\tau) 
                    - \frac\lambda\tau (F(\tau)\psi-F(\tau)\psi_\tau)
\]
Therefore
\[ (1+\frac\lambda\tau)(\psi-\psi_\tau) - (\phi_\tau-\phi)
   = \frac\lambda\tau (F(\tau)\psi-F(\tau)\psi_\tau)
\]
and taking norms and using the fact that $F$ is contraction we have
\[ \left| (1+\frac\lambda\tau)\|\psi-\psi_\tau\| -
   \|\phi_\tau-\phi\|\right|
   \leq 
   \| (1+\frac\lambda\tau)(\psi-\psi_\tau) - (\phi_\tau-\phi) \|
   \leq\frac\lambda\tau\|\psi-\psi_\tau\|
   \]
In particular
\[ (1+\frac\lambda\tau)\|\psi-\psi_\tau\| -
   \|\phi_\tau-\phi\| \leq\frac\lambda\tau\|\psi-\psi_\tau\|
\]
and therefore \(\|\psi-\psi_\tau\|\leq \|\phi-\phi_\tau\|\).

Now we estimate \(\|\phi-\phi_\tau\|\) in the simple case of a forward
Euler scheme. Note that 
\[ \phi-\phi_\tau = \lambda A \psi - \frac\lambda\tau(\Id-F(\tau))\psi
\]
and thus \(\|\phi-\phi_\tau\|\) measures a sort of consistency
error. For a single forward Euler step, \(F=F_E\) where \(F_E\) is
defined in \eqref{Brezis:onestep}. Thus 
\be
  \|\phi-\phi_\tau\| =\lambda \left\| A \psi -\frac1\st
  (\Id-S(\st))p(\psi)\right\|
  \rightarrow 0
\ee
as $\tau\rightarrow 0$ since $\frac{\Id-S(\st)}\st
p(\psi)\rightarrow Lp(\psi)=A\psi$.

The more general case of a $\nu$-stages Runge-Kutta scheme can be
carried out by induction following the procedure already applied
in the proofs of the previous lemmas.

We now use Theorem 3.2 of \cite{BP72} which, specialized to our case,
can be written as follows. Assume that $F(\tau):\mathrm{L}^1\rightarrow\mathrm{L}^1$
for $\tau>0$ is a family of contractions. Assume further that an
m-accretive operator $A$ is given and let $S(t)$ be the semigroup
generated by $A$. Assume further that the family $F(\tau)$ and the
operator $A$ are linked by the following formula 
\be
	\psi_{\tau}=\left(I+\frac{\lambda}{\tau}(I-F(\tau))\right)^{-1}\phi\rightarrow \psi=\left(I+\lambda A\right)^{-1}\phi
\ee
for each $\phi\in\mathrm{L}^1$. Then
\[ \lim_{n\rightarrow\infty} F\left(\frac{t}n\right)^n \phi = S(t)\phi 
   \qquad 
   \forall \phi\in\mathrm{L}^1. 
\]
\end{proof}

\section{Fully discrete relaxed scheme}                  
\label{section:fullydiscrete}
\setcounter{equation}{0}
\setcounter{figure}{0}
\setcounter{table}{0}

In order to complete the description of the scheme, we need to specify
the space discretization. We will use discretizations based on finite
differences, in order to avoid cell coupling due to the source terms.

Note that the IMEX technique reduces the integration to a cascade of
relaxation and transport steps. The former are the implicit parts of
\eqref{eq:stage:1} and \eqref{eq:stage:i}, while the transport steps
appear in the evaluation of the explicit terms $B^{(i)}$ in
\eqref{eq:stage:expl}. Since \eqref{eq:stage:1} and \eqref{eq:stage:i}
involve only local operations, the main task of the space discretization
is the evaluation of $\mathrm{div}(f)$, where we will exploit the
linearity of $f$ in its arguments.

\subsection{One dimensional scheme} \label{sec:oned:scheme}

Let us introduce a uniform grid on $[a,b]\subset\mathbb{R}$,
\(x_j=a-\frac{h}2+jh\) for $j=1,\ldots,n$, where $h=(b-a)/n$ is the
grid spacing and $n$ the number of cells. The fully discrete scheme may be written as
\be \label{eq:fd:imex}
z_j^{n+1} = z_j^n - \Delta t \sum_{i=1}^{\nu}
     \tilde{b}_i  \left(F^{(i)}_{j+1/2} - F^{(i)}_{j-1/2}\right)
    + \frac{\Delta t}{\varepsilon} \sum_{i=1}^{\nu} b_i
     g(z_j^{(i)}),  
\ee
where $F^{(i)}_{j+1/2}$ are the numerical fluxes, which are the only
item that we still need to specify. For convergence it is necessary to
write the scheme in conservation form. Thus, following
\cite{OS89}, we introduce the function $\hat{F}$ such that
\[
  f(z(x,t))=\frac1h \int_{x-h/2}^{x+h/2} \hat{F}(s,t) \mathrm{d}s
\quad\Rightarrow
  \frac{\partial f}{\partial x}(z(x_j,t)) = \frac1h
  \left(\hat{F}(x_{j+1/2},t)-\hat{F}(x_{j-1/2},t)\right).
\]
The numerical flux function $F_{j+1/2}$ approximates
$\hat{F}(x_{j+1/2})$. 

In order to compute the numerical fluxes, for each stage value, we
reconstruct boundary extrapolated data $z^{(i)\pm}_{j+1/2}$ with a
non-oscillatory interpolation method from the point values $z^{(i)}_j$ of the
variables at the center of the cells. Next we apply a monotone
numerical flux to these boundary extrapolated data.

To minimize numerical viscosity we choose the Godunov flux, which in
the present case of a linear system of equations reduces to the upwind
flux. In order to select the upwind direction we write the system in
characteristic form. The characteristic variables relative to the
eigenvalues $\varphi,-\varphi,0$ (in one space dimension $\varphi$ reduces to a scalar parameter) are respectively
\be \label{eq:char:var}
      U=\frac{v+\varphi w}{2\varphi}     \qquad      
      V=\frac{\varphi{}w-v}{2\varphi}  \qquad   
      W=u-w.
\ee
Note that $u=U+V+W$. Therefore the numerical flux in characteristic
variables is
\(F_{j+1/2}=(\varphi U^-_{j+1/2},-\varphi V^+_{j+1/2},0)\).

The accuracy of the scheme depends on the accuracy of the
reconstruction of the boundary extrapolated data. For a first order
scheme we use a piecewise constant reconstruction such that
\(U^-_{j+1/2}=U_j\) and \(V^+_{j+1/2}=V_{j+1}\). For higher order
schemes, we use ENO or WENO reconstructions of appropriate accuracy
(\cite{Shu97}).
 
For $\varepsilon\rightarrow0$ we obtain the relaxed scheme. Recall
from equation \eqref{relaxedRK}
that the
relaxation steps reduce to
\be \label{eq:relaxed:relaxation}
 w_j^{(i)}=p(u_j^{(i)}), \qquad v_j^{(i)}=-D\widehat\nabla w_j^{(i)},
\ee
where $\widehat\nabla$ is a suitable approximation of the one-dimensional
gradient operator. Thus the transport steps need to be applied only to
$u^{(i)}$
\be \label{eq:relaxed:transport}
 u_j^{(i)}=u_j^n -\lambda \sum_{k=1}^{i-1}\tilde{a}_{i,k}
   \left[
   \varphi \left( U^{(k)-}_{j+1/2}-U^{(k)-}_{j-1/2}\right)
   -\varphi \left( V^{(k)+}_{j+1/2}-V^{(k)+}_{j-1/2}\right)
   \right]
\ee
Finally, taking the last stage value and going back to conservative
variables, 
\be \label{eq:relaxed:RK}
\ba{ll} 
u_j^{n+1}=u_j^n 
    -\frac\lambda2
     \sum_{i=1}^\nu \tilde{b}_i 
         &\left([v^{(i)-}_{j+1/2}+v^{(i)+}_{j+1/2}-(v^{(i)-}_{j-1/2}+v^{(i)+}_{j-1/2})]\right.\\
              &\left.+\varphi[w^{(i)-}_{j+1/2}-w^{(i)+}_{j+1/2}-(w^{(i)-}_{j-1/2}-w^{(i)+}_{j-1/2})]
         \right)
\ea
\ee

We wish to emphasize that the scheme reduces to the time
advancement of the single variable $u$. Although the scheme is based
on a system of three equations, the construction is used only to select
the correct upwinding for the fluxes of the relaxed scheme and the
computational cost of each time step remains moderate.

\subsection{Non linear stability for the first order scheme}
The relaxed scheme in the first order case reduces to:
\be \label{eq:nonlin:stab}
  u_j^{n+1} = u_j^n 
        + \frac\lambda2 \left(\partial_x p(u^n)|_{j+1} - \partial_x p(u^n)|_{j-1}\right)
        + \frac\lambda2 \varphi \left( p(u_{j+1}^n)-2p(u_j^n)+p(u_{j-1}^n) \right)
\ee
We wish to compute the restrictions on $\lambda$ and $\varphi$ so that
the scheme is total variation non-increasing. We select the centered finite difference formula to
approximate the partial derivatives of $p(u)$; we drop the index $n$ and write $p_j$ for
$p(u^n_j)$. Define $\Delta_{j+1/2} = \frac{p_{j+1}-p_j}{u_{j+1}-u_{j}}$ and
observe that these quantities are always nonnegative since $p$ is
nondecreasing. 
We obtain
\be
\ba{lll}
  \TV(u^{n+1}) &= \sum_j &|u_j^{n+1}-u_{j-1}^{n+1}| = \\
       &\leq \sum_j 
       &\left\{ \frac\lambda{4h}\Delta_{j+3/2}|u_{j+2}-u_{j+1}|
               +\frac\lambda2\varphi\Delta_{j+1/2}|u_{j+1}-u_{j}|\right.\\
               &&\left.
               +\left(1-\lambda(\frac1{2h}+\varphi)\Delta_{j-1/2}\right)|u_{j}-u_{j-1}| \right.\\
               &&\left.+\frac\lambda2\varphi\Delta_{j-3/2}|u_{j-1}-u_{j-2}|
               +\frac\lambda{4h}\Delta_{j-5/2}|u_{j-2}-u_{j-3}|
       \right\}
\ea
\ee
provided that
\be \label{prestab:cond}
    (1-\lambda(\frac1{2h}+\varphi)\Delta_{j-1/2}) \geq 0 \qquad \forall j
\ee
Assuming that the data have compact support, we can
rescale all sums and finally get $\TV(u^{n+1})\leq\TV(u^{n}) $.
Taking into account the Lipschitz condition on $p$, the scheme is total
variation stable provided that 
\eqref{prestab:cond} is satisfied, i.e. that
\be \label{stab:cond}
    \Delta t
    \leq \frac{2h^2}{\mu}
         \frac1{1+2h\varphi} 
    \simeq \frac{(2-\delta)}{\mu} h^2 
\ee
where $\delta$ vanishes as $h$ does. We point out that the stability
condition is of parabolic type.
Finally, we observe that using one-sided approximations for the
partial derivatives of $p$ in the scheme \eqref{eq:nonlin:stab}, one
gets a stability condition involving the relation $\varphi>1/h$. This
would reintroduce in the scheme the constraint due to the stiffness in
the convective term that prompted the introduction of $\varphi$ in
\eqref{3eq}.
\subsection{Linear stability}
We study the linear stability of the schemes based on equations
\eqref{eq:relaxed:relaxation}, \eqref{eq:relaxed:transport} and
\eqref{eq:relaxed:RK} in the case when $p(u)=u$, by von Neumann
analysis. We substitute the discrete Fourier modes
$u_j^n=\rho^{n}e^{i(jk/N)}$ into the scheme, where $k$ is the wave
number and $N$ the number of cells. Let $\xi = k/N$,  we compute the amplification
factor $Z(\xi)$ such that $u_j^{n+1}=Z(\xi)u_j^n$. We can
consider $\xi$ as a continuous variable, since the
amplification factors for various 
choices of $N$ all lie on the curves obtained considering the variable
$\xi\in[0,2\pi]$.

\begin{figure}
\includegraphics[height=4.5cm]{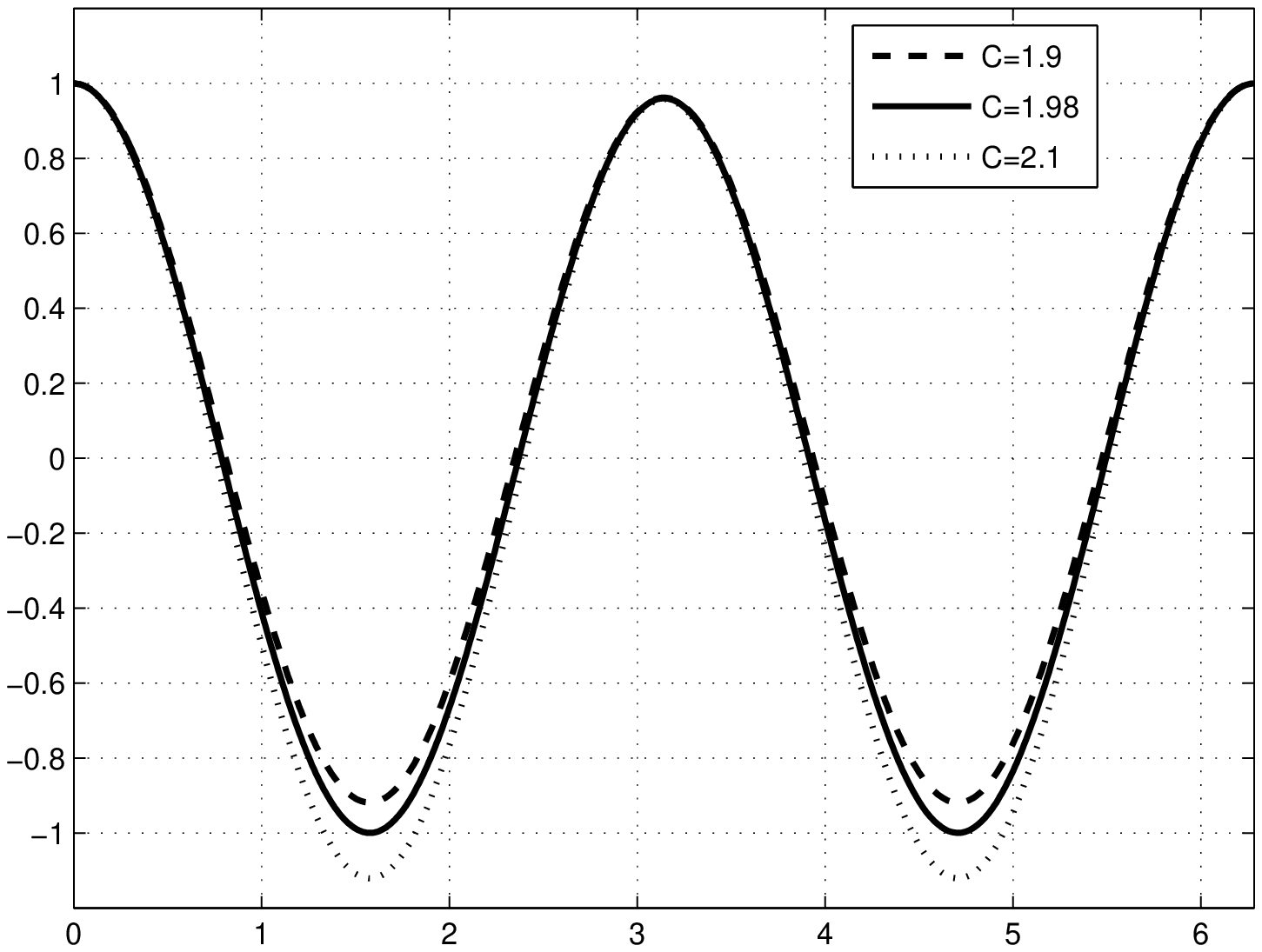}
\includegraphics[height=4.5cm]{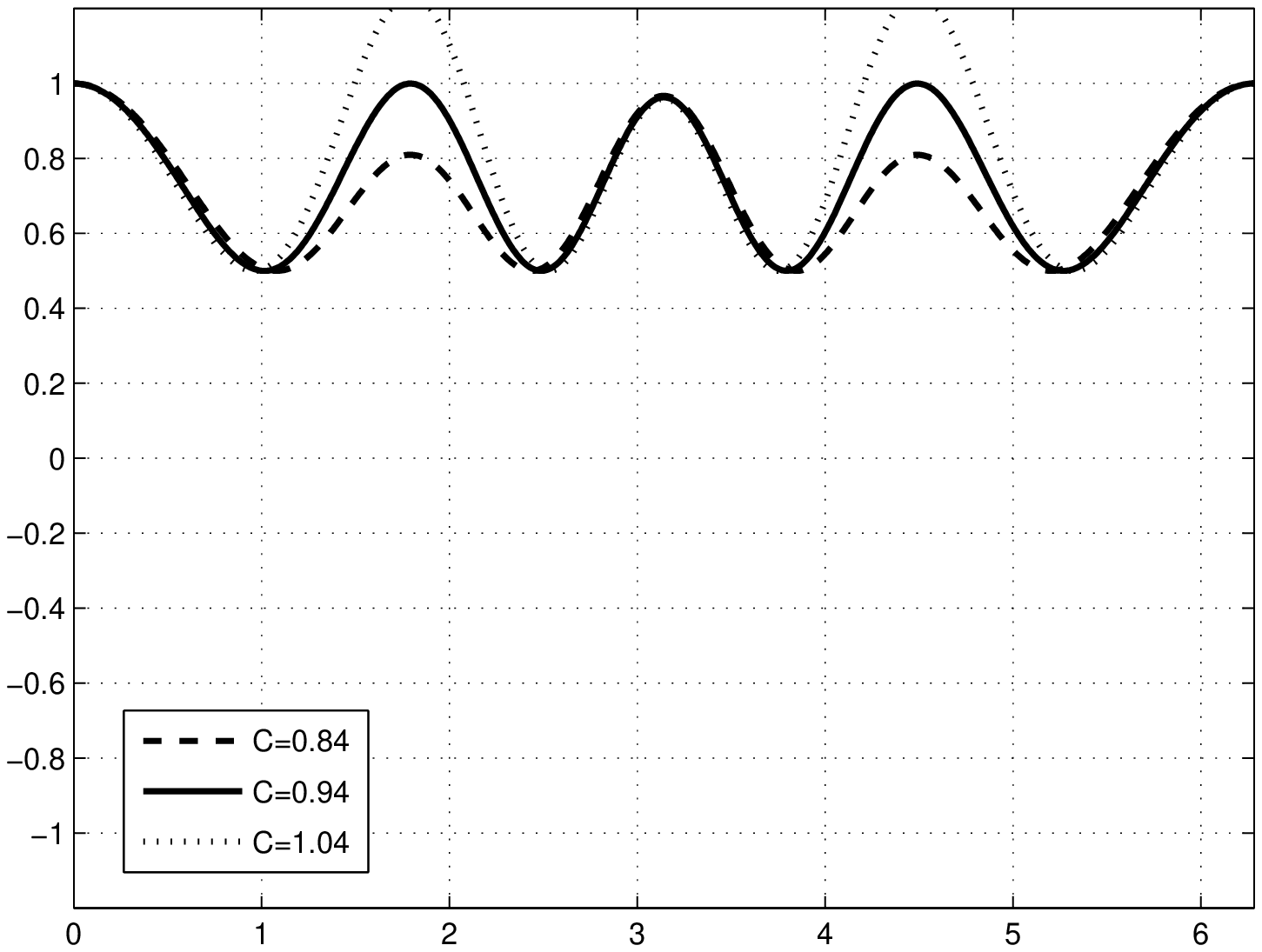}
\caption{Amplification factor for upwind spatial reconstruction 
coupled with forward Euler (left) and for upwind second order spatial
reconstruction  coupled with second order time integration (right).}
\label{fig:Up1Rk1}
\end{figure}

First we consider the same scheme studied in the previous section, for
comparison purposes. Using piecewise constant reconstructions in space
and forward Euler time integration, the amplification factor is
$Z(\xi)=1+M(\xi)$, where
\[M(\xi)=\frac\lambda{h}\left(\cos(\xi)-1\right)\left(\cos(\xi)+1+h\varphi\right).\]
 \(M(\xi)\) takes maximum value $0$ and attains its minimum at the point
$\xi^*$ such that $\cos(\xi^*)=-\varphi{}h/2$.
Stability requires that \(M(\xi^*)\geq-2\), i.e.
\[ 1+\frac\lambda{h}\left(\frac{\varphi^2h^2}{4}-1\right)-\lambda\varphi\left(\frac{\varphi{h}}{2}+1\right) \geq-1\]
and recalling that $\lambda=\mathrm{\Delta t}/h$,
\be  \label{eq:stab:lin:U1R1}
  \Delta t 
  \leq \frac{2h^2}{\left(1+\frac{\varphi{h}}{2}\right)^2}
  \simeq 2\left(1-\varphi{h}\right) h^2
\ee
This gives a CFL condition of the form \(\Delta t \leq 2(1-\delta)
h^2\) where $\delta=O(h\varphi)$ (see figure \ref{fig:Up1Rk1}).
These results are in very good agreement with those of
the nonlinear analysis performed in the previous section.

Now we consider higher order spatial reconstructions coupled with
forward Euler timestepping.  $M$ takes the form
\[ M(\xi,\gamma) = \frac\lambda{h} \left[f_1(\cos(\xi)) + \gamma f_2(\cos(\xi))\right] \]
where $\gamma=h\varphi$. Since $\gamma$ is small, we compute the
critical points $\xi^*$ of $M(\xi,0)$. For stability we thus require
that $-2\leq{}M(\xi^*,\gamma)\leq0$.

We consider a piecewise linear and a WENO reconstruction. The first
one is computed along characteristic variables using the upwind slope,
while the gradient of $p(u)$ is computed with centered
differences. The WENO reconstruction is fifth order accurate and is
obtained by setting to $1$ the smoothness indicators and the gradient
of $p(u)$ is computed with the fourth order centered difference
formula.

For the piecewise linear reconstruction, we have
that
\[ M(\xi) = -\frac{\lambda}{h}
            \left[ (\cos^2(\xi)-1)(\cos(\xi)-2) + h\varphi(\cos(\xi)-1)^2 
            \right].
\]
and therefore
\[ \Delta t 
   \leq \frac{2h^2}{\frac{20+14\sqrt{7}}{27} + \frac{8+2\sqrt{7}}{9}\varphi{}h} 
   \simeq 0.94(1 - 1.44\varphi{}h) h^2
\]
For the WENO reconstruction $M(\xi,\gamma)$ can be easily computed and
we get
\[ \Delta t \leq 0.79 (1-0.13\varphi{}h) h^2\]

Now we wish to extend our results to the case of higher order
Runge-Kutta schemes. Since both the equation and the scheme are
linear, the amplification factors for the Runge-Kutta schemes of order
$2$ and $3$ used here are respectively  
\[ \ba{l} 
    Z_{(2)}(\xi) = 1+ M(\xi) + \frac{M(\xi)^2}{2} \\
    Z_{(3)}(\xi) = 1+ M(\xi) + \frac{M(\xi)^2}{2} + \frac{M(\xi)^3}{6}.
\ea\]
where $M(\xi)$ is the function appearing in the amplification factor
relevant to the chosen spatial reconstruction. We have that
\[ \ba{l}
      Z_{(2)}'(\xi) = M'(\xi) (1+M(\xi)) \\
      Z_{(3)}'(\xi) = M'(\xi) (1+M(\xi)+\frac{M(\xi)^2}{2})
\ea\]
and therefore 
the critical points are 
the points $\xi^*$ such that
$M'(\xi^*)=0$.

In the Runge-Kutta $2$ case the
stability constraint $\|Z_{(2)}(\xi^*)\|\leq1$ reduces to the CFL
condition for 
the forward Euler  scheme.
For Runge-Kutta 3,
\( \|Z_{(3)}(\xi^*)\| \leq 1 \), 
provided that
\[ M(\xi^*)\geq\tilde{s}\simeq-2.51 \]
Notice that this is less restrictive than the Euler and RK2 schemes
for which the stability requirement is $M(\xi^*)\geq-2$.

For the Runge-Kutta 3 scheme with linearized WENO of order 5,
we have
\[ \Delta t 
   \leq \frac{-\tilde{s}h^2}{2.51 + 0.33\varphi{}h} 
   \simeq (1-.1325\varphi{}h) h^2
\]

\begin{figure}
\includegraphics[height=4.5cm]{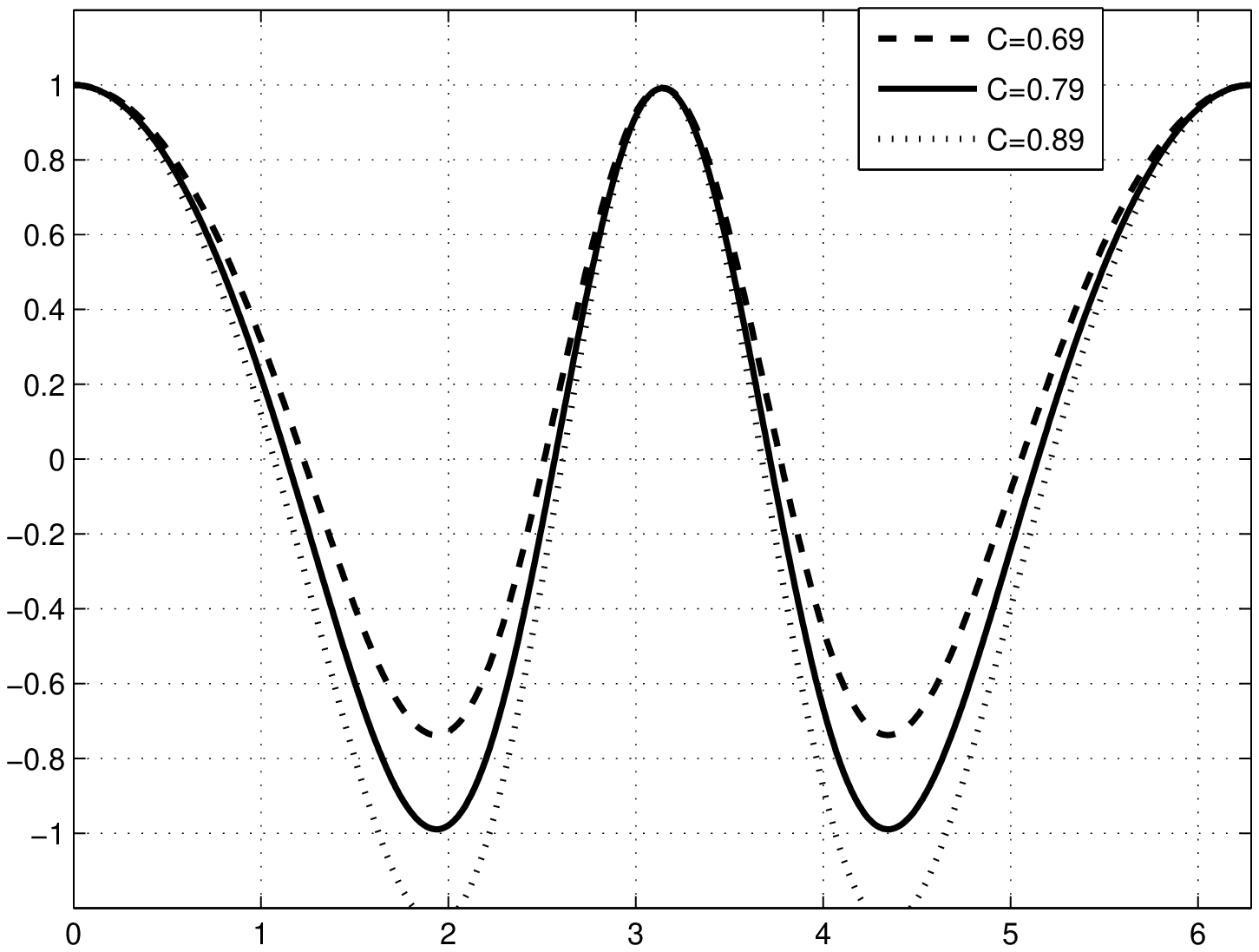}
\includegraphics[height=4.5cm]{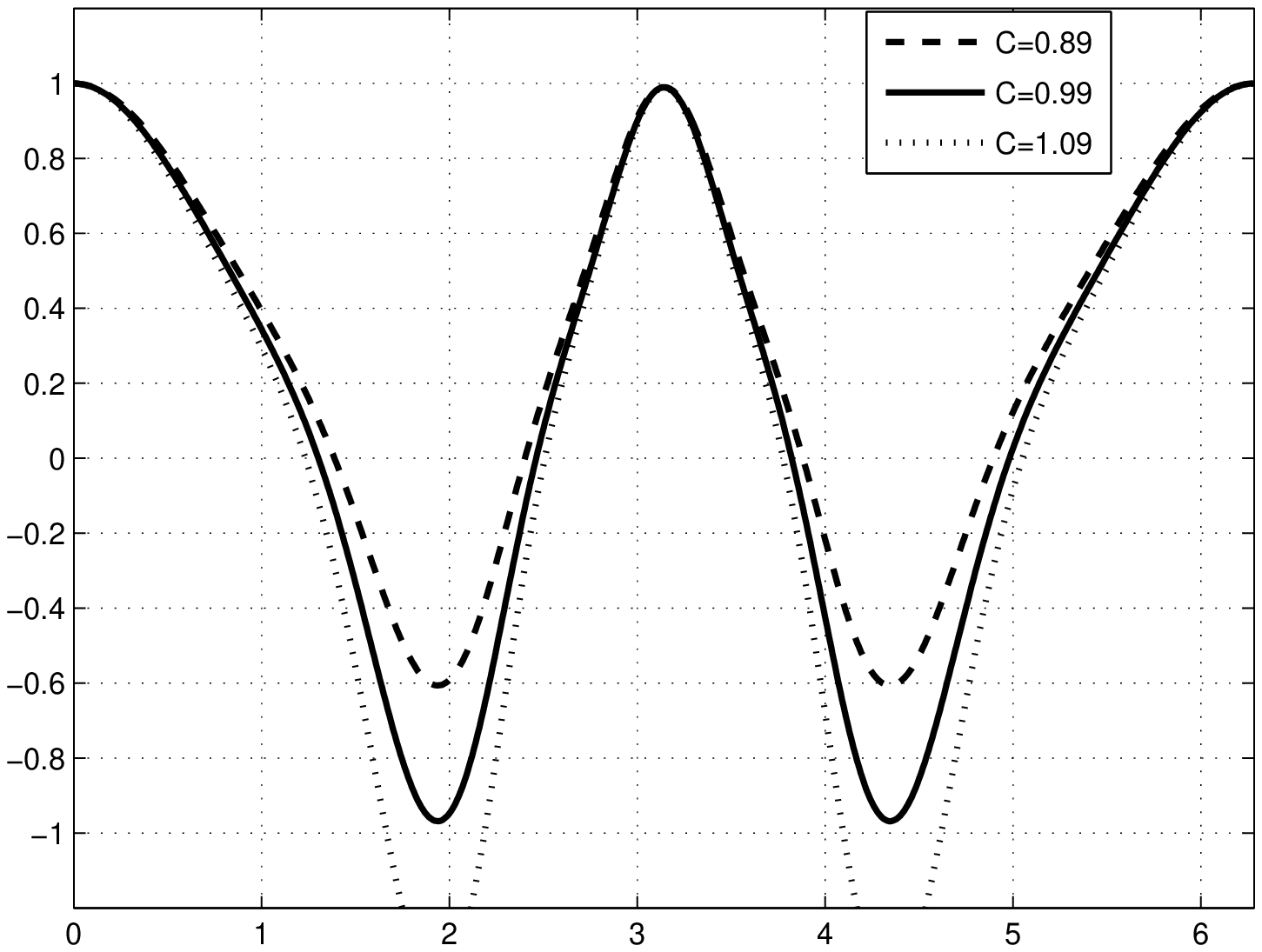}
\caption{Amplification factors $Z$ for WENO reconstructions of order $5$
coupled with first order (left) and third order (right) time integration.}
\label{fig:WE3RK13}
\end{figure}

Table \ref{table:RK} summarizes the stability results obtained in this
section listing the values of the constant $C$ that appear in the
stability restriction $\Delta t\leq C(1-C_1\varphi{h})h^2$. Figures
\ref{fig:Up1Rk1} and \ref{fig:WE3RK13} contain the amplification
factors $Z(\xi)$ for $\varphi=1$ and $h=10^{-2}$ for various choices
of spatial reconstructions and time integration schemes. Each of them
contains the curve corresponding to the value of $C$ reported in Table
\ref{table:RK} and two other close-by values.

\begin{table}[htb]\scriptsize
\hfil
\begin{tabular}{|c|c|c|c|}
\hline
                & RK1 & RK2 & RK3\\
\hline
P-wise constant & 2   & 2   & 2.51\\
\hline
P-wise linear   & 0.94&0.94 &     \\
\hline
WENO5           & 0.79&0.79 & 1   \\
\hline
\end{tabular}
\hfil
\caption{}\label{table:RK}
\end{table}

\subsection{Boundary conditions} \label{sec:bc}

Different boundary conditions can be implemented. Here we describe how
to implement Neumann boundary conditions, considering for simplicity
the one-dimensional case.

 We first add $g$ ghost points on each side
of the computational domain $[a,b]$, where $g$ depends on the order of the
spatial reconstruction.
We find a polynomial $q(x)$ of degree $d$ passing through the points
$(x_i,u_i)$ for $i=1,\ldots,d$ and having prescribed derivative at the
boundary point $x_{1/2}=a$. (The degree $d$ is determined by the
accuracy of the scheme that one wants to obtain and should match the
degree of the reconstruction procedures used to obtain $U^{\pm}_j$ and
$V^{\pm}_j$.)  This polynomial is then used to set the values
$u_{-i}=q(x_{-i})$ of the ghost points for $i=0,1,g-1$. One operates
similarly at the right edge of the computational domain.

We also used periodic boundary conditions, which can be implemented
with obvious choice of the values $u_i$ at the ghost points.

\subsection{Multi-dimensional scheme} \label{sec:highd:scheme}

An appropriate numerical approximation of \eqref{3eq} in
$\mathbb{R}^d$ that generalizes the scheme described in Section
\ref{sec:oned:scheme} can be obtained by additive dimensional
splitting. We consider the relaxed scheme, i.e. $\varepsilon=0$ and
for the sake of simplicity, let us focus on the square domain
\([a,b]\times[a,b]\subset\mathbb{R}^2\). 
Here we shall
describe the generalization of the scheme defined by equations
\eqref{eq:relaxed:relaxation}, \eqref{eq:char:var},
\eqref{eq:relaxed:transport} and \eqref{eq:relaxed:RK} to the case of
two space dimensions.

Without loss of generality, we consider a uniform grid in \([a,b]\times[a,b]\subset\mathbb{R}^2\) such that
\(\vec{x}_{i,j}=(x_i,y_j)=(a-h/2,a-h/2)+i(h,0)+j(0,h)\) for
\(i,j=1,2,\ldots,n\) and $h=(b-a)/n$. 

In the present case, $u$ and $w$ are one-dimensional variables, while
$\vec{v}=(v_{(1)},v_{(2)})$ is now a field in $\mathbb{R}^2$. First we
observe that the relaxation steps
\eqref{eq:relaxed:relaxation} are easily generalized for $d>1$. For
the transport steps, one has to evolve in time the system 
\begin{equation}
\frac{\partial}{\partial t}
\left(\ba{c} u\\v_{(1)}\\v_{(2)}\\w\ea\right)
+
\frac{\partial}{\partial x}
\left[\ba{cccc}
0 & 1 & 0 & 0 \\
0 & 0 & 0 & \varphi^2 \\
0 & 0 & 0 & 0 \\
0 & 1 & 0 & 0
\ea\right]
\left(\ba{c} u\\v_{(1)}\\v_{(2)}\\w\ea\right)
+
\frac{\partial}{\partial y}
\left[\ba{cccc}
0 & 0 & 1 & 0 \\
0 & 0 & 0 & 0 \\
0 & 0 & 0 & \varphi^2 \\
0 & 0 & 1 & 0
\ea\right]
\left(\ba{c}
 u\\v_{(1)}\\v_{(2)}\\w
\ea\right)
=0
\end{equation}
The semidiscretization in space of the above equation can be written as
\[ \frac{\partial z_{i,j}}{\partial t} = 
    -\frac1h \left(F_{i+1/2,j}-F_{i-1/2,j}\right)
    -\frac1h \left(G_{i,j+1/2}-G_{i,j-1/2}\right),
\]
where $F$ and $G$ are the numerical fluxes in the $x$ and $y$
direction respectively and can be written as
\[
F_{i+1/2,j} = F(z^+_{i+1/2,j},z^-_{i+1/2,j})
\qquad
G_{i,j+1/2} = G(z^+_{i,j+1/2},z^-_{i,j+1/2})
\]
The fluxes in the two directions are computed separately. We
illustrate the computation of the flux $F$ along the $x$ direction. 
We note that only the
field $v_{(1)}$ appears in the differential operator along this
direction. The third component of the flux is zero and thus we have
three independent characteristic variables, namely
\[ U_{(1)}=\frac{\varphi{}w+v_1}{2\varphi}  \qquad  
   V_{(1)}=\frac{\varphi{}w-v_1}{2\varphi}  \qquad   
   W=u-w,
\]
which correspond respectively to the eigenvalues
$\varphi,-\varphi,0$. At this point the numerical fluxes can be easily
evaluated by upwinding. We proceed similarly for the numerical flux
$G$ that depends on the characteristic variables $U_{(2)},V_{(2)},W$.

Denote by $U^\pm_{i+1/2,j}$ the reconstructions of \(U_{(1)}(\cdot,y_j)\)
at the point $(x_i+h/2,y_j)$. This involves a reconstruction of the
restriction of $U_{(1)}$ to the line $y=y_i$ and can be obtained with any of
the one-dimensional techniques mentioned in Section \ref{sec:oned:scheme}. 
Similarly, denote $U^\pm_{i,j+1/2}$ the reconstructions of \(U_{(2)}(x_i,\cdot)\)
at the point $(x_i,y_j+h/2)$.
Now, formulas \eqref{eq:relaxed:transport} and \eqref{eq:relaxed:RK} become respectively
\be \label{eq:relaxed:transport:2d}
\ba{ll}
 u_{i,j}^{(l)}=u_{i,j}^n -\lambda \sum_{m=1}^{l-1}\tilde{a}_{l,m}
   &\left[
   \varphi \left( U^{(m)-}_{i+1/2,j}-U^{(m)-}_{i-1/2,j}\right)
   -\varphi \left( V^{(m)+}_{i+1/2,j}-V^{(m)+}_{i-1/2,j}\right) \right.\\
   &\left.
   \varphi \left( U^{(m)-}_{i,j+1/2}-U^{(m)-}_{i,j-1/2}\right)
   -\varphi \left( V^{(m)+}_{i,j+1/2}-V^{(m)+}_{i,j-1/2}\right)
   \right]
\ea
\ee
and
\be \label{eq:relaxed:RK:2d}
\ba{ll}
 u_{i,j}^{n+1}=u_{i,j}^n -\lambda \sum_{l=1}^{\nu} 
   \varphi\tilde{b}_l
   &\left[
   \left( U^{(l)-}_{i+1/2,j}-V^{(l)+}_{i+1/2,j} \right)
   - \left(U^{(l)-}_{i-1/2,j} -V^{(l)+}_{i-1/2,j}\right)\right.\\
   &\left.
   \left( U^{(l)-}_{i,j+1/2}-V^{(l)+}_{i,j+1/2} \right)
   - \left(U^{(l)-}_{i,j-1/2} -V^{(l)+}_{i,j-1/2}\right)
   \right]
\ea
\ee

The generalization to $d>2$ and rectangular domains is now trivial. We
stress once again that no two-dimensional reconstruction is used, but
only $d$ one-dimensional reconstructions are needed. Finally,
boundary conditions can be implemented direction-wise with the
same techniques used in the one-dimensional case.

\section{Numerical results}                  
\label{section:numerical}
\setcounter{equation}{0}
\setcounter{figure}{0}
\setcounter{table}{0}

We performed several numerical tests of our relaxed schemes. First we
tested convergence for a linear diffusion equation 
with periodic and Neumann boundary conditions for initial data giving
rise to smooth solutions. 
Next, numerical tests were also performed on the porous media equation
$u_t=(u^m)_{xx}$, $m=2,3$, both in one and two dimensions.

\subsection{Linear diffusion}
\label{sec:numer:lin:per}

For the first test we considered the linear problem
\[
\left\{
\begin{array}{lll}
\displaystyle
\frac{\partial u}{\partial t}(x,t) 
      = \frac{\partial^2 u}{\partial x^2} u(x,t) &x\in[0,1]\\
\displaystyle
u(x,0)=u_0(x)                                     &x\in[0,1]
\end{array}
\right.
\] 
First we used periodic boundary conditions with $u_0(x)=\cos(2\pi x)$,
so that $u(x,t)=\cos(2\pi x)e^{-4\pi^2t}$ is an exact solution. Then
we used Neumann boundary conditions $u_x(0)=u_x(1)=1$ with initial
data $u_0(x)=x+cos(2\pi x)$, so that $u(x,t)=x+\cos(2\pi x)e^{-4\pi^2t}$
is an exact solution.

\begin{table}
\begin{centering}
\footnotesize
\begin{tabular}{|c|l|l|l|l|l|}
\hline
 & N=40 & N=80 & N=160 & N=320 & N=640 \\
\hline
ENO2, RK1 &  2.012e-03 &  5.6378e-04 &  1.0736e-04 &  1.5539e-05 & 2.5065e-06 \\
\hline
ENO3, RK2 &  1.9066e-06 &  2.3057e-07 &  5.6115e-08 &  8.6904e-09 &  1.1905e-09\\
\hline
ENO4, RK2 & 7.7517e-06 &  5.7082e-07 &  3.3507e-08 &  1.4978e-09 &  7.0725e-11\\
\hline
ENO5, RK3 & 1.3864e-08 &  6.0259e-10 &  2.2121e-11 &  7.4454e-13 & 2.3803e-14 \\
\hline
ENO6, RK3 & 1.5538e-08 & 8.5661e-10 & 1.446e-11 & 1.7111e-13 & 1.5311e-15\\
\hline
WENO3, RK2 & 1.9799e-03  &  5.1278e0-4  &  1.4332e-04  &  2.1488e-05  &   7.512e-08\\
\hline
WENO5, RK3 & 1.5892e-07 &  4.8069e-09 &    1.59e-10 &  5.2337e-12  & 1.6758e-13\\
\hline
\end{tabular}

\medskip
\begin{tabular}{|c|l|l|l|l|l|}
\hline
 & N=40 & N=80 & N=160 & N=320 & N=640 \\
\hline
ENO2, RK1   & 1.3973      & 1.8354      & 2.3926      & 2.7886      & 2.6322\\
\hline
ENO3, RK2   & 5.9501      & 3.0477      & 2.0388      & 2.6909      & 2.8678\\
\hline
ENO4, RK2   & 3.8987      & 3.7634      & 4.0905      & 4.4836      & 4.4045\\
\hline
ENO5, RK3  & 6.8124       & 4.524      & 4.7677      & 4.8929      & 4.9671\\
\hline
ENO6, RK3  & 5.9907       & 4.181      & 5.8885       & 6.401      & 6.8043\\
\hline
WENO3, RK2 & 0.56648       & 1.949      & 1.8391      & 2.7376      & 8.1601\\
\hline
WENO5, RK3  & 2.9595      & 5.0471       & 4.918       & 4.925      & 4.9649\\
\hline
\end{tabular}
\caption{$\mathrm{L}^1$ norms of the error and convergence rates for the linear diffusion equation
with periodic boundary conditions, with smooth initial data} 
\label{tab:errori:lin:per}
\end{centering}
\end{table}

\begin{table}
\begin{centering}
\footnotesize
\begin{tabular}{|c|l|l|l|l|l|}
\hline
 & N=40 & N=80 & N=160 & N=320 & N=640 \\
\hline
ENO2, RK1 &  2.1965e-03 &  5.7152e-04 &  1.4301e-04  &  2.32e-05  &  4.743e-06\\
\hline
ENO3, RK2 &  2.0621e-06 &  2.2641e-07 &  6.7935e-08 &  8.8255e-09 &  1.2339e-09\\
\hline
ENO4, RK2 & 8.1764e-06 &  5.4431e-07 &  3.6974e-08 &  1.3686e-09  &  8.335e-11\\
\hline
ENO5, RK3 & 1.5484e-07 &  4.4163e-09 &  1.2405e-10 &  3.7803e-12 &  1.1669e-13\\
\hline
WENO3, RK2 & 1.9092e-03 &  4.4225e-04 &  1.2914e-04 &  4.5037e-06 &  7.4526e-08\\
\hline
WENO5, RK3 & 2.5048e-07 &  4.9279e-09 &  1.4776e-10 &  4.7482e-12 &  1.4948e-13\\
\hline
\end{tabular}

\medskip
\begin{tabular}{|c|l|l|l|l|l|}
\hline
 & N=40 & N=80 & N=160 & N=320 & N=640 \\
\hline
ENO2, RK1   & 1.4361   &     1.9424    &    1.9987    &     2.624     &   2.2902\\
\hline
ENO3, RK2   & 6.1004    &    3.1871    &    1.7367    &    2.9444    &    2.8385\\
\hline
ENO4, RK2   & 3.9763    &     3.909    &    3.8798    &    4.7558    &    4.0373\\
\hline
ENO5, RK3  & 5.6626   &     5.1317    &    5.1539     &   5.0362    &    5.0178\\
\hline
WENO3, RK2 & 1.2624    &      2.11    &    1.7759    &    4.8417   &     5.9172\\
\hline
WENO5, RK3  &   4.9122  &     5.6676  &     5.0597   &    4.9597 &   4.9893 \\
\hline
\end{tabular}
\caption{$\mathrm{L}^1$ norms of the error and convergence rates for the linear diffusion equation
with Neumann boundary conditions, with smooth initial data.}
\label{tab:errori:lin:neu}
\end{centering}
\end{table}

We tested the numerical schemes defined by equations
\eqref{eq:char:var}, \eqref{eq:relaxed:relaxation}, 
\eqref{eq:relaxed:transport} and \eqref{eq:relaxed:RK} with various
degrees of accuracy for the spatial reconstructions and time-stepping
operators.
We used ENO spatial reconstructions of degrees from $2$ to $6$ and
WENO reconstructions of degrees $3$ and $5$. The time-stepping
procedures chosen are IMEX Runge-Kutta schemes of Section
\ref{section:semidiscrete} of accuracy chosen to
match the accuracy of the spatial reconstruction. Since stability forces
the parabolic restriction $\Delta t \leq C h^2$, an IMEX scheme of
order $m$ was coupled with a spatial ENO/WENO reconstruction of
accuracy $p$ such that $p\leq2m$, obtaining a scheme of order $p$.

We computed the numerical solution of the diffusion equation with
final time $t=0.05$ with $N=40, 80, 160, 320, 640$ grid points and
computed the $\mathrm{L}^1$ norm of the difference between the numerical and
the exact solution. The results are in Tables
\ref{tab:errori:lin:per} for the periodic boundary conditions and
\ref{tab:errori:lin:neu} for the Neumann boundary conditions. One can
see that the expected convergence 
rates are reached, even if the combination of the
WENO reconstruction of accuracy 3 and the IMEX scheme of second order
reach the predicted error reduction only on very fine grids.

\subsection{Porous media equation}
\label{sec:numer:porous}

On the porous media equation \eqref{eq:degparab} with \(p(u)=u^m\) 
we performed a test
proposed in \cite{GJ71}. We took \(m=2,3\) and initial data of
class $C^1$ as follows:
\be \label{eq:cosquadro}
 u(x,0)= \left\{ \ba{ll} \cos^2(\pi x/2) & |x|\leq1 \\ 0 & |x|>1 \ea
\right.
\ee
The computational domain is 
$\{|x|\leq3\}\subset\mathbb{R}$ and the boundary conditions are
periodic; the CFL constant is taken as $C=0.25$.

\begin{table}
\begin{centering}
\footnotesize
\begin{tabular}{|c|l|l|l|l|l|}
\hline
 & N=60 & N=180 & N=540 & N=1620 \\
\hline
ENO2, RK1 &  2.6365e-04 &  1.9898e-05 &   2.049e-06  &  2.076e-07 \\
\hline
ENO3, RK2 &  1.9605e-05 &  6.0423e-07 &  2.4141e-08  &  8.9729e-10 \\
\hline
ENO4, RK2 &  1.2127e-05 &  2.967e-07  &  9.9925e-09  &  3.5781e-10 \\
\hline 
ENO5, RK3 &  4.694e-06  &  1.719e-07  &  6.3248e-09  &  2.4447e-10 \\
\hline 
ENO6, RK3 &  4.1099e-06 &  1.4711e-07 &  5.3992e-09  &  2.0849e-10 \\
\hline
WENO3, RK2 & 1.5871e-04 &  1.0448e-05 &  4.3463e-07  & 8.8767e-09\\
\hline
WENO5, RK3 & 7.5662e-06 &  4.6049e-07 &  7.4746e-09  & 2.7985e-10\\
\hline
\end{tabular}

\medskip
\begin{tabular}{|c|l|l|l|l|l|}
\hline
 & N=60 & N=180 & N=540 & N=1620\\
\hline
ENO2, RK1   & 2.8243 &       2.352  &     2.0692  &      2.084 \\
\hline
ENO3, RK2   & 5.1899 &      3.1672  &      2.931  &     2.9968 \\
\hline
ENO4, RK2   & 5.6271 &      3.3774  &     3.0865  &     3.0307 \\
\hline
ENO5, RK3  &  6.491 &      3.0103   &     3.006   &    2.9611 \\
\hline
ENO6, RK3  &  6.612 &      3.0311   &    3.0083   &     2.962\\
\hline
WENO3, RK2 &  3.2863 &      2.4765  &     2.8942  &     3.5418  \\
\hline
WENO5, RK3  & 6.0565 &      2.5479  &     3.7509  &     2.9902  \\
\hline
\end{tabular}
\caption{$\mathrm{L}^1$ norms of the error and convergence rates for
  the porous media equation periodic boundary conditions, with initial  
data of class $C^1$.} 
\label{tab:errori:J2}
\end{centering}
\end{table}

Since the initial data has compact support and is Lipschitz continuous,
the solution will be of compact support for every $t\geq0$, but
will develop a discontinuity in $u_x$ at some finite time $\tau>0$
(see \cite{Aro70}).

As was shown in \cite{Aro70}, the solution with
the initial condition we chose has a front that does not move for
$t<0.034$. We therefore chose a final time of the simulation
$t_{\mbox{\tiny fin}}=0.03$ to prevent the formation of the
singularity of $u_x$ from affecting the order of convergence.
We used as reference solution the one obtained numerically with
$N=4860$ grid points and computed the $\mathrm{L}^1$ norms of the errors of the
solutions with $N=60,180,540,1620$ grid points. The results are
presented in Table \ref{tab:errori:J2}.

First of all one verifies that the degree of regularity of the solution
poses a limit on the order of convergence of the schemes: therefore the
schemes we tested perform at best as third order schemes, as
confirmed by the data in Table \ref{tab:errori:J2}.
Still, high order  schemes yield smaller error on a given grid. This can
be of practical importance in problems where one does not have the
freedom of choosing the
number of grid points, as in digital image analysis, where non-linear
degenerate diffusion equations are sometimes used as filters
for contour enhancement (see \cite{BV04}). 

\begin{figure}
\hfil
\includegraphics[height=4.3cm]{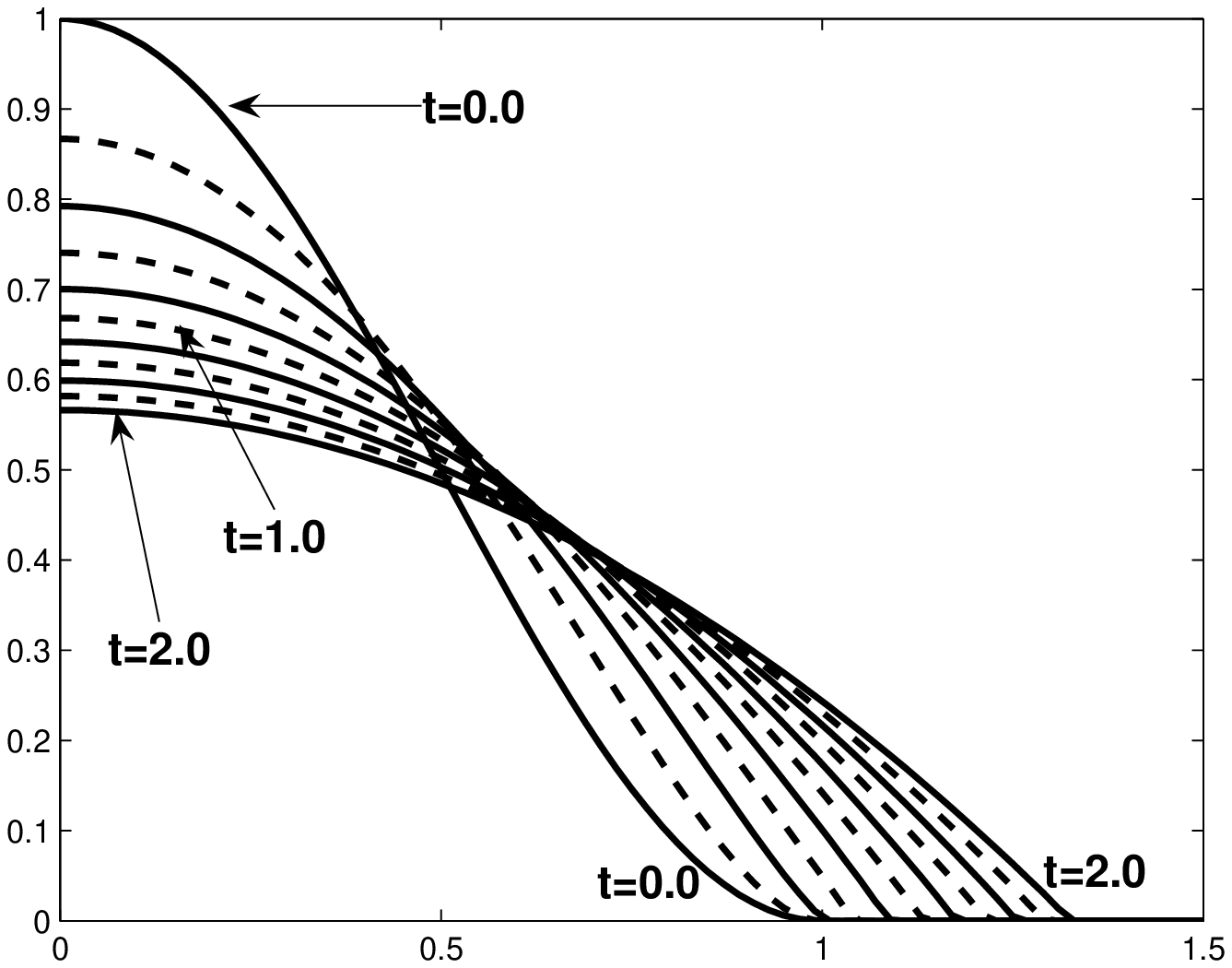}
\includegraphics[height=4.3cm]{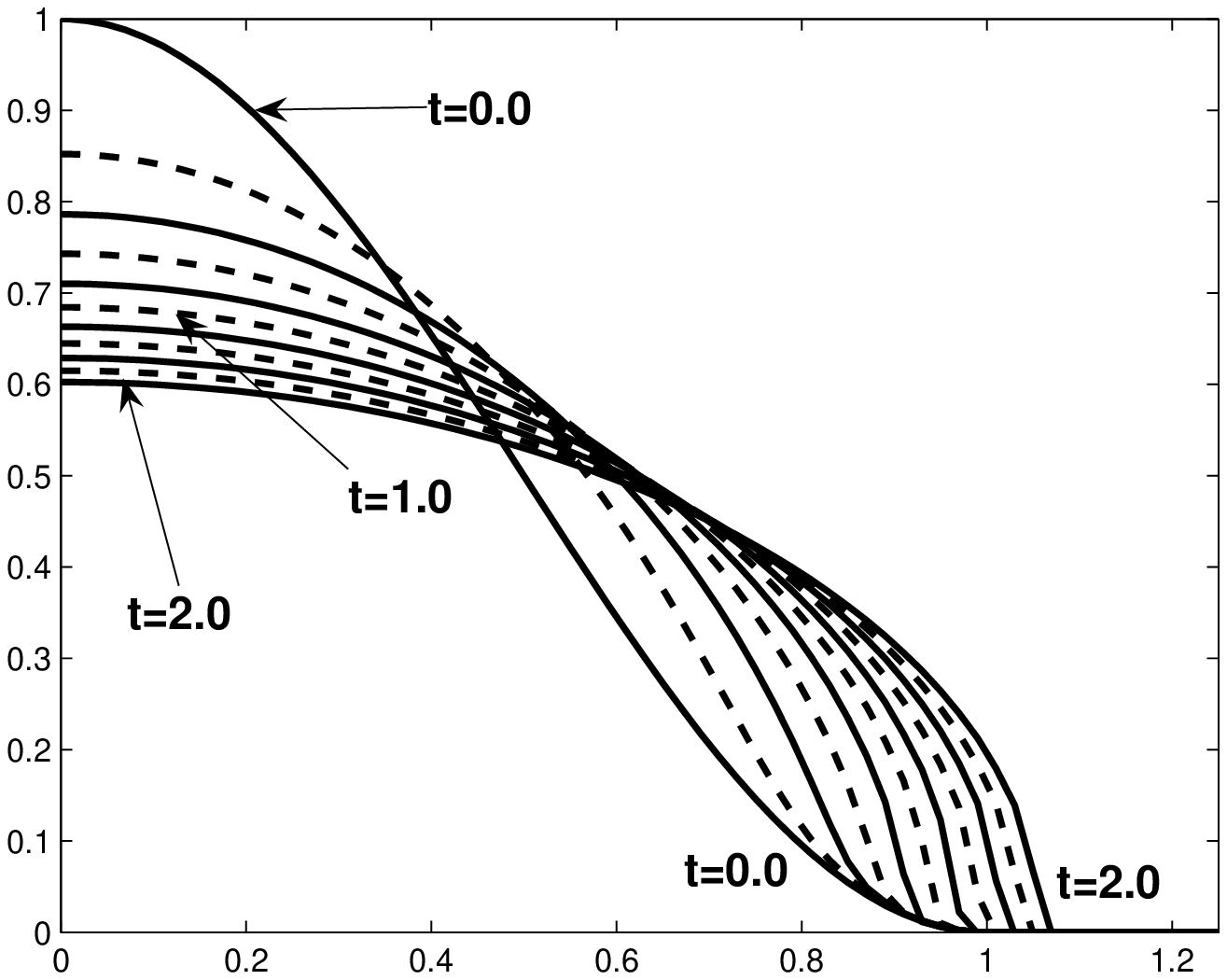}
\hfil
\caption{Snapshots of the numerical solutions for the porous media
  equation with $p(u)=u^2$ (left) and $p(u)=u^3$ (right). Initial data
  are chosen according to (\ref{eq:cosquadro}) and the numerical
  solutions are represented at times $t=0,0.2,\ldots,2.0$. The
  solutions are obtained with the spatial WENO reconstruction of order
  $5$ and the RK3 time integrator.}
\label{fig:ucubo}
\end{figure}
In Figure \ref{fig:ucubo} we show the numerical solution for the
porous media equation with $p(u)=u^2$ and $p(u)=u^3$, with the initial
data (\ref{eq:cosquadro}) and $t\in[0,2]$. It can be appreciated that
a front (i.e. a discontinuity of $\frac{\partial{u}}{\partial{}x}$)
develops at a finite time and then it travels at finite speed.

\begin{figure}
\begin{centering}
\epsfig{file=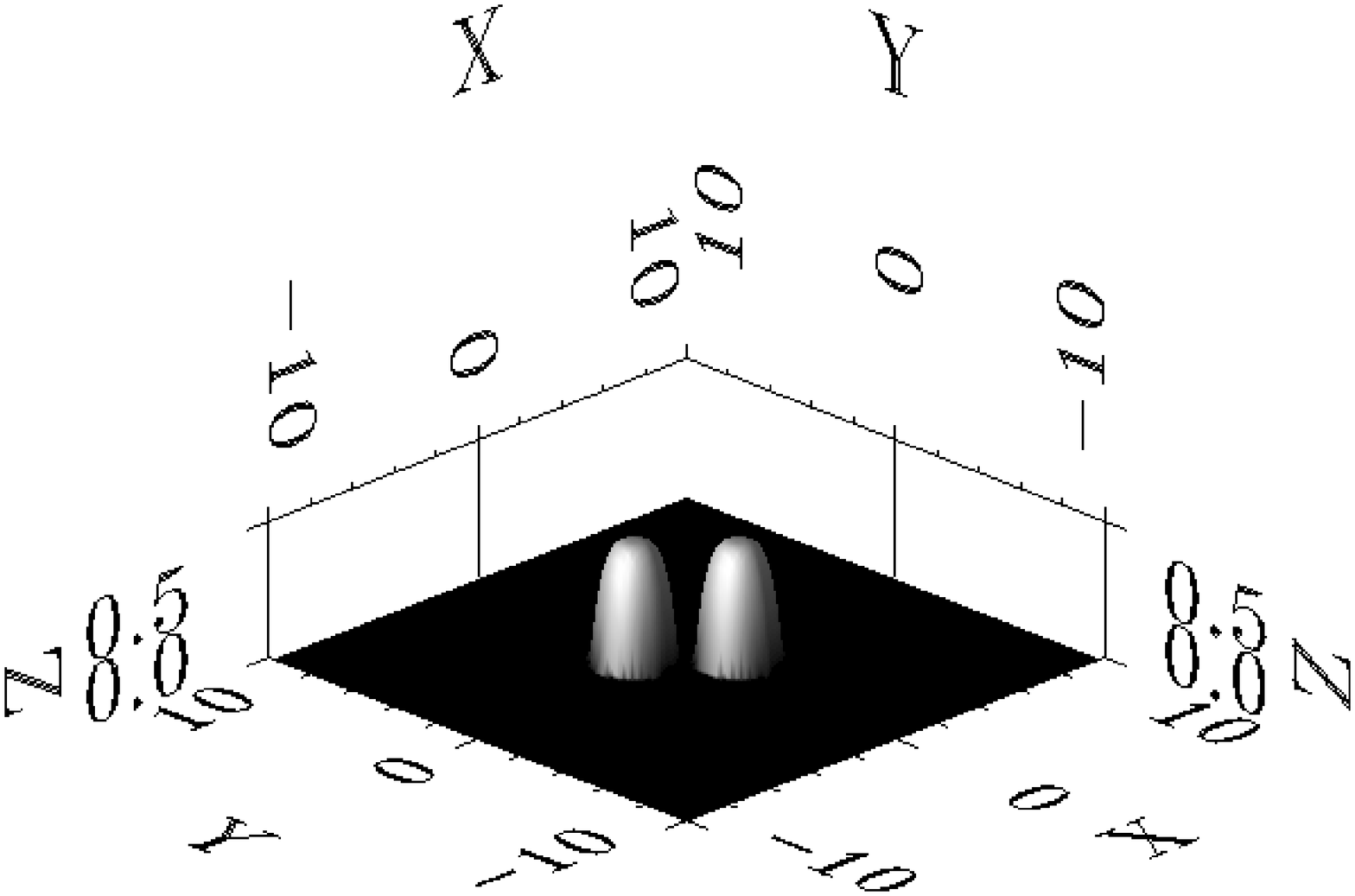,height=.45\textwidth,angle=-90}\hfill
\epsfig{file=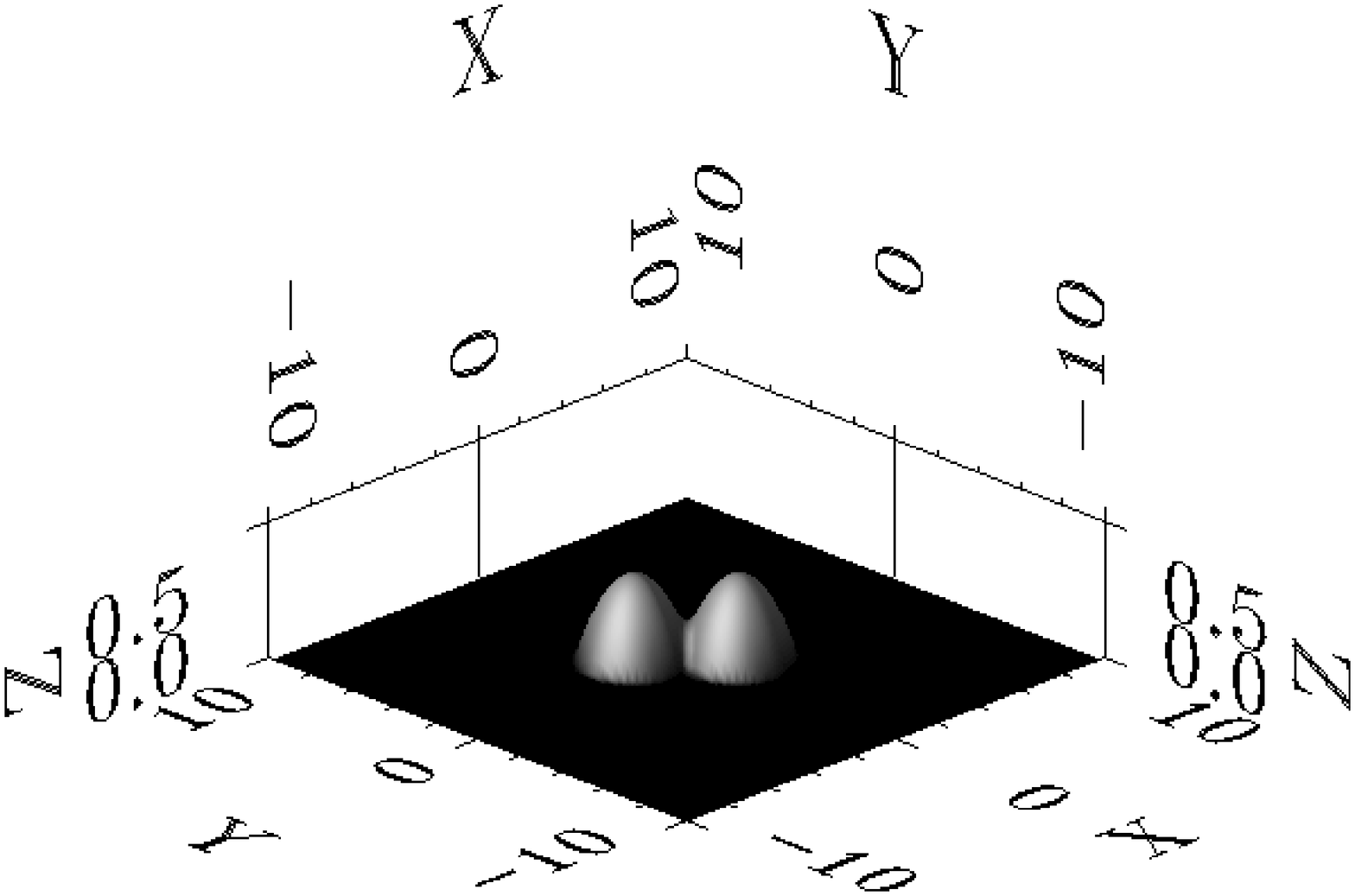,height=.45\textwidth,angle=-90}\\
\epsfig{file=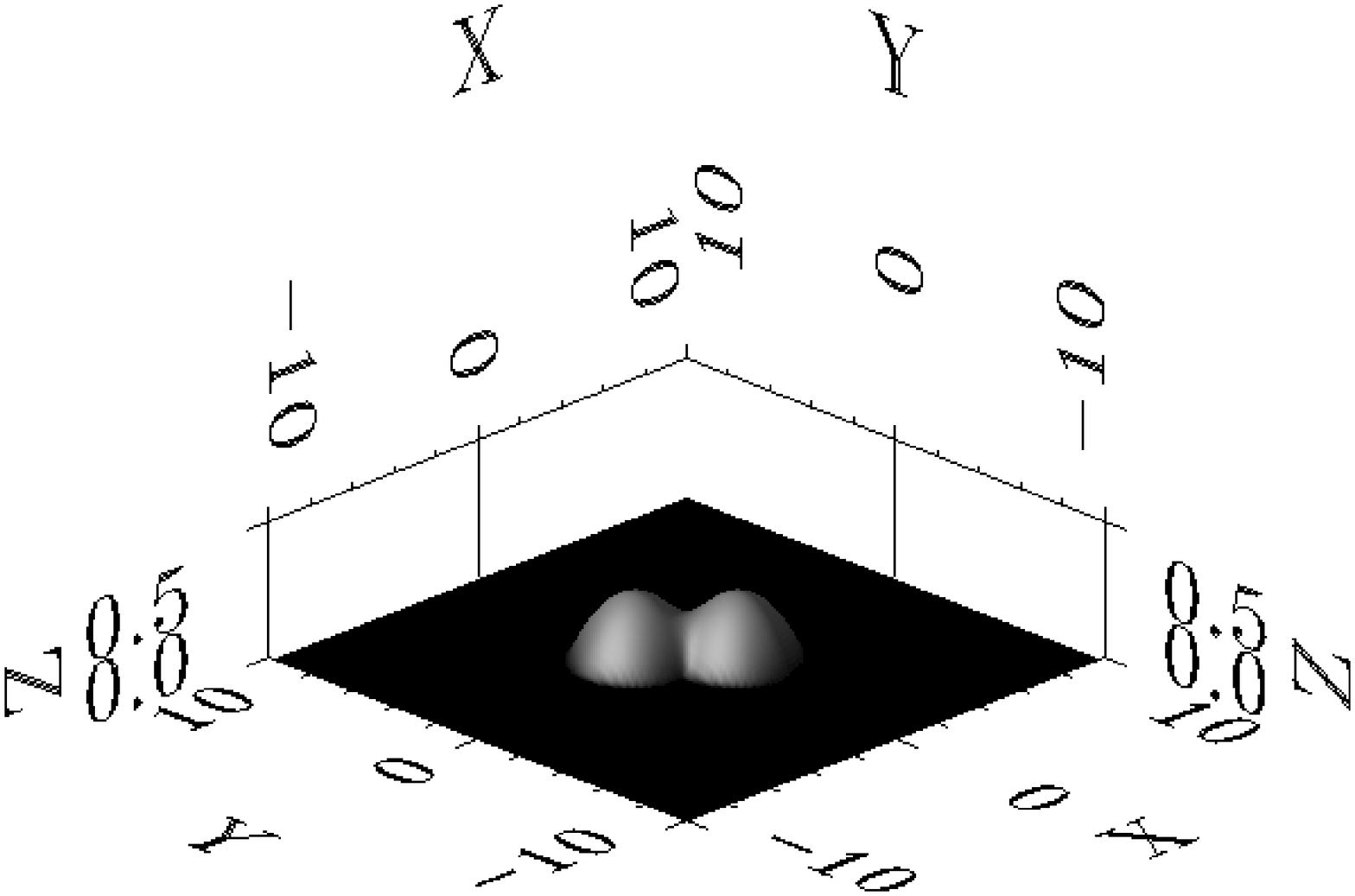,height=.45\textwidth,angle=-90}\hfill
\epsfig{file=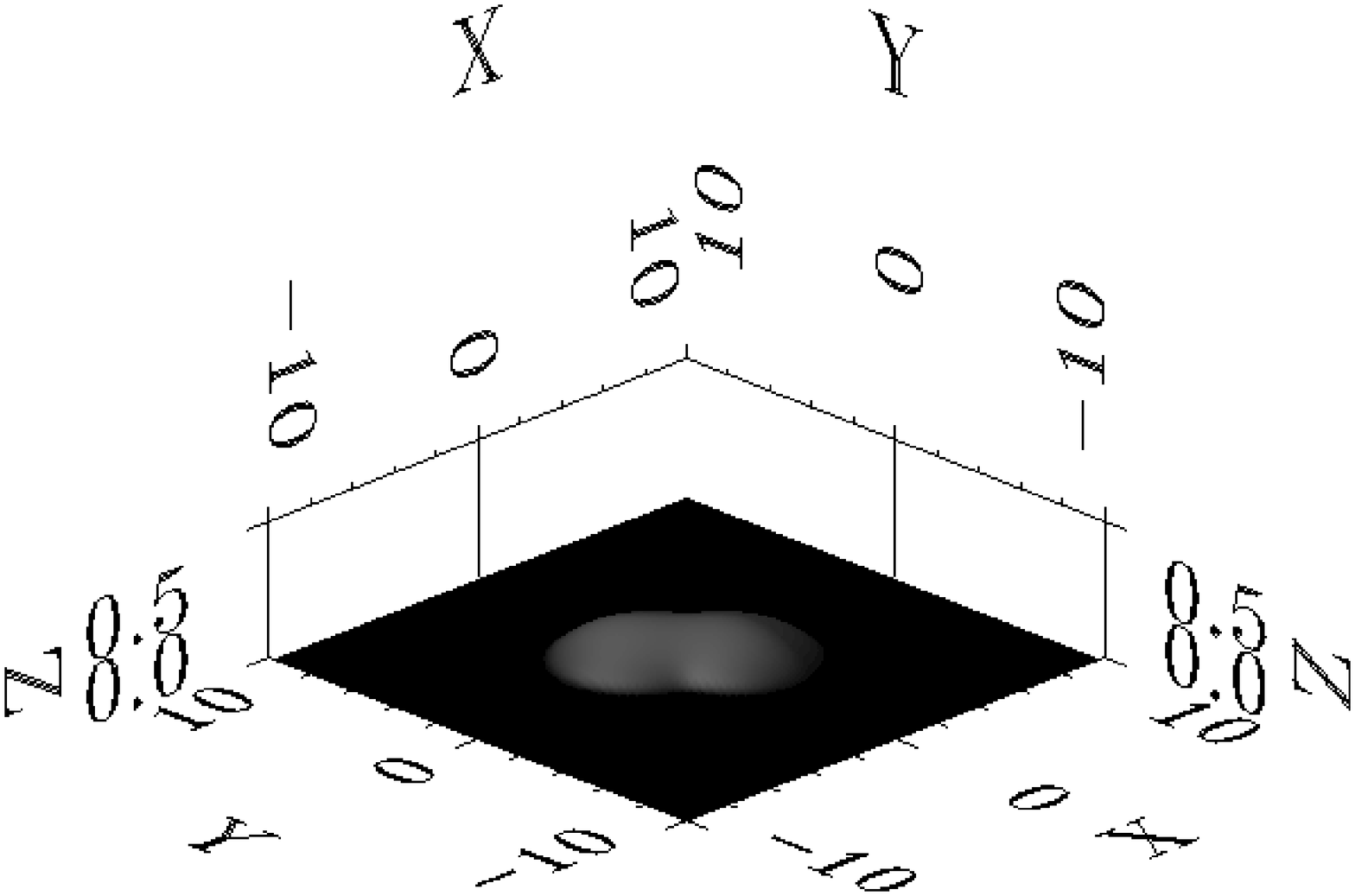,height=.45\textwidth,angle=-90}\\
\end{centering}
\caption{The numerical solution of the porous media equation on a
square regular grid with compactly supported initial data. From top
left to bottom right, we show the
numerical solution at times $t=0,0.5,1.0,4.0$.}
\label{fig:2d}
\end{figure}

We present a numerical simulation for the two-dimensional porous media
equation \eqref{eq:degparab} with $p(u)=u^2$. 
 We chose an initial data $u_0(x,y)$ given by two bumps
with periodic boundary conditions on $[-10,10]\times[-10,10]$. The
large domain ensures that the compact support of the solution is
still contained in the computational domain at the final time of the
calculation. The numerical approximation at different time levels is
shown in Figure \ref{fig:2d}. 

We can note that the symmetries of the initial data are preserved and
the solution seems to be unaffected by the dimensional splitting of
the two-dimensional scheme.

\section{Conclusions}                  
\label{section:conlusions}
\setcounter{equation}{0}
\setcounter{figure}{0}
\setcounter{table}{0}

We have proposed and analyzed relaxed schemes for nonlinear degenerate
parabolic equations.

By using suitable discretization in space and time, namely ENO/WENO
non-oscillatory reconstructions for numerical fluxes and IMEX
Runge-Kutta schemes for time integration, we have obtained a class of
high order schemes. We have developed a theoretical convergence
analysis for  the semidiscrete scheme; furthermore we studied
stability for the fully discrete schemes. Our computational results
suggest that our schemes converge with the predicted
rate.

Finally, we point out that these schemes can be easily implemented on
parallel computers. Some preliminary results and details are reported
in \cite{enumath05}. In particular the schemes involve only linear
matrix-vector operations and the execution time scales linearly when
increasing the number of processors.

Our numerical approach can be easily extended
also to more general problems, as nonlinear convection-diffusion
equations or nonlinear parabolic systems. Some of these applications
will appear in a forthcoming paper.

\bibliographystyle{plain}
\bibliography{0604572}

\end{document}